\documentclass[12pt]{amsart}
\usepackage{tikz-cd}
\usepackage{comment}
\usepackage{latexsym,amssymb,amsmath}
\usepackage{hyperref}
\usepackage{comment}
\usepackage{mathrsfs}
\usepackage{amsmath,amscd}
\usepackage[enableskew]{youngtab}
\usepackage[all]{xy}
\usepackage{empheq}
\usepackage{amsmath}
\usepackage{empheq}
\usepackage{xcolor}
\usepackage{amssymb,hyperref}
\usepackage{url}
\usepackage{relsize}

\excludecomment{verlong}
\includecomment{vershort}
\excludecomment{noncompile}
\usepackage{tikz}
\usetikzlibrary{matrix}
\usepackage{combelow}

\newtheorem{theorem}{Theorem}[section]
\newtheorem{lemma}[theorem]{Lemma}
\newtheorem{proposition}[theorem]{Proposition}
\newtheorem{corollary}[theorem]{Corollary}

\newtheorem{rem}[theorem]{Remark}
\newtheorem{defn}[theorem]{Definition}
\newtheorem{exmp}{Example}[section]

\numberwithin{equation}{section}

\DeclareMathOperator{\id}{id}

\DeclareMathOperator{\ev}{ev}

\DeclareMathOperator{\rank}{{rank}}

\DeclareMathOperator{\rev}{rev}

\usepackage{tikz}
\usetikzlibrary{matrix}

\usepackage{combelow}

\usepackage{amsmath,amsthm,amsfonts,amssymb}
\usepackage{amsmath, amsthm, amssymb, amscd, amsfonts}
\usepackage[all]{xy}
\usepackage{amssymb, amsthm, amsmath}
\usepackage{amsfonts}

\usepackage{color}
\DeclareMathOperator{\ad}{ad}

\newcommand{\ra}{\rightarrow}

\newcommand{\ot}{\otimes}
\newcommand{\co}{\mathcal O}
\newcommand{\xra}{\xrightarrow}
\newcommand{\mtc}{\mathcal}
\newcommand{\cs}{\mtc S}

\newcommand{\lam}{\lambda}

\newcommand{\Lam}{\Lambda}

\newcommand{\al}{\alpha}
\newcommand{\eps}{\epsilon}

\newcommand{\ul}{\underline}

\newcommand{\lh}{\leftharpoonup}
\newcommand{\whb}{{\widehat{(H, \mtc B)}}}
\newcommand{\hb}{{(H, \mtc B)}}

\newcommand{\dw}{\downarrow}
\newcommand{\uw}{\uparrow}

\newcommand{\ch}{\chi}
\newcommand{\mtr}{\mathrm}

\newcommand{\ncm}{\newcommand}
\ncm{\np}{\newpage}
\ncm{\ebl}{\end{thebibliography}}
\ncm{\bbl}{\begin{thebibliography}}
\ncm{\chd}{_{ _{\ch}}}
\ncm{\ald}{_{ _{\al}}}
\newcommand{\blam}{\Lam}
\ncm{\cP}{\mathcal{P}}
\ncm{\ei}{e_i}
\ncm{\eij}{e_{i,\;j}}
\ncm{\bt}{\begin{theorem}}
\ncm{\bdef}{\begin{defn}}
\ncm{\edf}{\end{defn}}
\ncm{\et}{\end{theorem}}
\ncm{\bc}{\begin{corollary}}
\ncm{\bl}{\begin{lemma}}
\ncm{\el}{\end{lemma}}
\ncm{\bpf}{\begin{proof}}
\ncm{\epf}{\end{proof}}
\ncm{\ec}{\end{corollary}}
\ncm{\ord}{\mtr{ord}}
\ncm{\er}{\end{rem}}
\ncm{\br}{\begin{rem}}
\ncm{\bn}{\begin}

\ncm{\bp}{\begin{proposition}}
\ncm{\ep}{\end{proposition}}
\ncm{\bd}{
\begin{document}}
\ncm{\ed}{\end{document}}
\ncm{\beq}{\begin{equation}}
\ncm{\beqn}{\begin{equation*}}
\ncm{\eeq}{\end{equation}}
\ncm{\eeqn}{\end{equation*}}
\ncm{\bea}{\begin{eqnarray}}
\ncm{\eea}{\end{eqnarray}}
\ncm{\beanon}{\begin{eqnarray*}}
\ncm{\eeanon}{\end{eqnarray*}}\ncm{\ek}{\eps|_K}\ncm{\diez}{\#}
\ncm{\bwt}{\bowtie}
\ncm{\cC}{\mtc{C}}\ncm{\cc}{\mtc{C}}
\ncm{\cX}{\mtc{X}}
\ncm{\wt}{\widetilde}
\ncm{\sg}{\sigma}
\ncm{\Rep}{\mathrm{Rep}}
\ncm{\X}{\mathcal{X}}
\ncm{\cA}{\mathcal{A}}
\ncm{\HKer}{\mtr{HKer}}
\ncm{\LKER}{\mtr{LKer}}
\ncm{\aad}{\mtr{ad}}
\newcommand{\mbf}{\mathbb F}
\ncm{\Dr}{\mtr{D}}
\ncm{\cD}{{\mathcal{D}}}\ncm{\cd}{{\mathcal{D}}}\ncm{\ce}{{\mathcal{E}}}
\ncm{\G}{\mathcal{G}}
\ncm{\Dc}{\mtc{D}}
\ncm{\E}{\mtc{E}}
\ncm{\fp}{\mtr{FPdim}}
\ncm{\Vc}{\mtr{Vec}}
\ncm{\cK}{\mtc{K}}
\ncm{\cM}{\mtc{M}}
\ncm{\cE}{\mtc{E}}
\ncm{\cS}{\mtc{S}}
\newcommand{{\ipr}}{i'}
\newcommand{\tomega}{\widetilde{\omega}}
\ncm{\cop}{\mtr{cop}}
\ncm{\op}{\mtr{op}}
\ncm{\chr}{character }\ncm{\ck}{\mtc{K}}
\ncm{\bw}{\bwt}
\ncm{\hker}{\mtr{HKer}}
\ncm{\bx}{\boxtimes}
\ncm{\blue}{\textcolor[rgb]{.00, .00, 1.00}}
\ncm{\red}{\textcolor[rgb]{1.00, .00, .00}}
\ncm{\green}{\textcolor[rgb]{.50, 0.20, .90}}
\ncm{\bne}{\begin{enumerate}}
\ncm{\ene}{\end{enumerate}}
\ncm{\lker}{\mtr{LKer}}
\ncm{\md}{\medbreak}
\ncm{\rep}{\Rep}\ncm{\ind}{\mtr{ind}}
\ncm{\mdn}{\md\noindent}
\ncm{\dd}{$}
\ncm{\up}{^}
\newcommand{\tcs}{\text}
\newcommand{\mbb}{\mathbb B}
\newcommand{\vs}{\mathbb V}
\newcommand{\sth}{suppose that\;}
\newcommand\rad{\operatorname{rad}}
\newcommand{\itm}{\item}
\newcommand{\dbd}{$$}
\newcommand{\mol}{\mtr{mod}}
\newcommand{\ro}{\rho}
\newcommand{\irr}{\mathrm{Irr}}
\newcommand{\mbc}{\mathbb C}
\newcommand{\mbs}{\mathbb S}
\newcommand{\mbz}{\mathbb Z}
\newcommand{\ct}{\mtc T}
\newcommand{\sm}{\setminus}
\newcommand{\epl}{^{+}}
\newcommand{\sbsq}{\subseteq}
\newcommand{\sbs}{\subset}
\newcommand{\cco}{\mtr{co}}
\newcommand{\cz}{\mathcal{Z}}
\newcommand{\dual}{^{*}}
\newcommand{\Gm}{\Gamma}
\ncm{\cY}{\mtc{Y}}
\newcommand\ZZ{{\mathbb Z}} 
\newcommand{\bab}{\color{DarkOrchid}{}}
\newcommand{\eab}{\normalcolor{}}
\newcommand{\subs}{\subsection}
\newcommand{\cv}{\mtc{V}}
\newcommand{\grn}{\green}
\newcommand{\dt}{\delta}
\newcommand{\ccf}{\mathrm{ {CF}(\cc)}}
\newcommand{\cce}{\mathrm{ {CE}(\cc)}}
\newcommand{\cecc}{\mathrm{ {CE}(\cc)}}
\newcommand{\cecd}{\mathrm{ {CE}(\cd)}}
\newcommand{\kk}{\Bbbk}
\newcommand{\otL}{\ot_{L}}
\newcommand{\otl}{\ot_{L}}
\newcommand{\unpsi}{1_{\psi}}
\newcommand{\epsi}{e_{\psi}}
\newcommand{\ephi}{e_{\phi}}
\newcommand{\ech}{e_{\ch}}
\newcommand{\nleftcid}{\text{left normal  coideal subalgebra}}
\newcommand{\dimL}{\dim_{\kk}L}
\newcommand{\cl}{\mtc L}
\newcommand{\mj}{\mtc J}
\newcommand{\tl}{\tilde L}
\newcommand{\tL}{\tilde L}
\newcommand{\tpsi}{\tilde(\psi)}
\newcommand{\tmx}{\tilde{\mtc X}}
\newcommand{\zlh}{\mathrm{ZL}}
\newcommand{\ba}{\mathrm A}
\newcommand{\bv}{\mathrm V}
\newcommand{\zhopf}{\mtc{Z}_{\mtr{Hopf}}}
\newcommand{\lstar}{L^{*}}
\newcommand{\ldstar}{L^{**}}
\newcommand{\mstar}{M^{*}}
\newcommand{\mdstar}{M^{**}}
\newcommand{\lkera}{\lker_{A}}
\newcommand{\mdprime}{M''}
\newcommand{\ldprime}{L''}
\newcommand{\cm}{\mtc M}
\newcommand{\ccm}{\mathcal M}
\newcommand{\cn}{\mathcal N}
\newcommand{\ccn}{\mathcal N}
\newcommand{\rx}{\mtr{Rex}}
\newcommand{\cca}{\ca}
\newcommand{\ih}{\underline{\mtr{Hom}}}
\newcommand{\cih}{\underline{\mtr{coHom}}}
\newcommand{\hm}{\mtr{ {Hom}}}
\newcommand{\cov}{\mtr{coev}}
\newcommand{\rora}{\rho^{\mtr{ra}}}
\newcommand{\rola}{\rho^{\mtr{la}}}
\newcommand{\cx}{\mtc X}
\newcommand{\cZ}{\cz}
\newcommand{\ca}{\cA}
\newcommand{\stat}{\noindent}
\newcommand{\bfa}{{\bf A}}
\newcommand{\unu}{\mathbf{1}}
\newcommand{\barzu}{{\bar {  Z}(\unu)}}
\newcommand{\idx}{\id_X}
\newcommand{\lprime}{L'}
\newcommand{\mprime}{M'}
\newcommand{\nat}{ \mtr{{  Nat}}}
\newcommand{\ft}{\mtc F_\lam}
\newcommand{\rhau}{\rightharpoonup}
\newcommand{\lhau}{\leftharpoonup}
\newcommand{\cf}{\mathrm{ {CF}}}
\newcommand{\cfc}{\mathrm{{CF}}(\cc)}
\newcommand{\csu}{\overline{\mathfrak{  C}}}
\newcommand{\cfcc}{{\mathrm{CF}(\cc)}}
\newcommand{\catfcc}{\mathrm{ {CF}}(\cc)}
\newcommand{\cfcd}{{\mathrm{CF}(\cd)}}
\newcommand{\cfd}{{\mathrm{CF}(\cd)}}
\newcommand{\czcc}{{\cz(\cc)}}
\newcommand{\czcd}{{\cz(\cd)}}
\newcommand{\czt}{{\cz(\cz(\cc))}}
\newcommand{\enx}{\mtr{  End}}
\newcommand{\runu}{R(\unu)}
\newcommand{\bdfn}{\bn{defn}}
\newcommand{\edfn}{\end{defn}}
\newcommand{\deltax}{\delta_X}
\newcommand{\deltav}{\delta_V}
\newcommand{\repcca}{\rep_\cc(A)}
\newcommand{\xotay}{X \ot_A Y}
\newcommand{\xoty}{X \ot Y}
\newcommand{\votw}{V \ot W}
\newcommand{\votaw}{V \ot_A W}
\newcommand{\dimax}{\dim_AX}
\newcommand{\dimccx}{\dim_\cc(X)}
\newcommand{\dimcca}{\dim_\cc(A)}
\newcommand{\dimccv}{\dim_\cc(V)}
\newcommand{\dima}{\dim_A}
\newcommand{\biga}{A}
\newcommand{\comp}{\mathbb C}
\newcommand{\tehtaa}{\theta_A}
\newcommand{\tetaa}{\theta_A}
\newcommand{\ida}{\id_A}
\newcommand{\hma}{\hm_A}
\newcommand{\hmcc}{\hm_\cc}
\newcommand{\fv}{F(V)}
\newcommand{\fw}{F(W)}
\newcommand{\ota}{\ot_A}
\newcommand{\repza}{\rep_\cc^0(A)}
\newcommand{\epsa}{\eps_A}
\newcommand{\bndefn}{\bn{defn}}
\newcommand{\edefn}{\end{defn}}
\newcommand{\bdefn}{\bn{defn}}
\newcommand{\vld}{V^{*}}
\newcommand{\vldd}{V^{**}}
\newcommand{\xld}{X^{*}}
\newcommand{\xldd}{X^{**}}
\newcommand{\yld}{Y^{*}}
\newcommand{\yldd}{Y^{**}}
\newcommand{\aldu}{A^{*}}
\newcommand{\aldd}{A^{**}}
\newcommand{\ia}{\mtr{i}_A}
\newcommand{\aota}{A\ot A}
\newcommand{\idv}{\id_V}
\newcommand{\ld}{^*}
\newcommand{\repg}{\rep(G)}
\newcommand{\thetav}{\theta_V}
\newcommand{\tta}{\theta_A}
\newcommand{\muv}{\mu_V}
\newcommand{\muw}{\mu_W}
\newcommand{\dimcc}{\dim(\cc)}
\newcommand{\chii}{\chi_i}
\newcommand{\chistar}{\ch_{i^*}}
\newcommand{\chj}{\ch_j}
\newcommand{\chm}{\ch_m}
\newcommand{\chn}{\ch_n}
\newcommand{\dimvi}{\dim(V_i)}
\newcommand{\mtcd}{Q}
\newcommand{\mtca}{\mtc A}
\newcommand{\lamcd}{\lam_\cd}
\newcommand{\fpdimcd}{\fp(\cd)}
\newcommand{\laml}{\lam_L}
\newcommand{\apm}{A//M}
\newcommand{\apl}{A//L}
\newcommand{\repapm}{\rep(\apm)}
\newcommand{\repapl}{\rep(\apl)}
\newcommand{\dimvj}{\dim(V_j)}
\newcommand{\dvi}{\dim(V_i)}
\newcommand{\dvj}{\dim(V_j)}
\newcommand{\sumjtom}{\sum_{j=0}^m}
\newcommand{\sumitom}{\sum_{i=0}^m}
\newcommand{\sij}{s_{ij}}
\newcommand{\sji}{s_{ji}}
\newcommand{\dxj}{d_j}
\newcommand{\dxi}{\di }
\newcommand{\dimka}{\dim_{\kk}(A)}
\newcommand{\dimk}{\dim_{\kk}}
\newcommand{\blaml}{\blam_L}
\newcommand{\sumjtor}{\sum_{j=0}^r}
\newcommand{\dimkl}{\dim_{\kk}(L)}
\newcommand{\mtcjl}{\mtc J_L}
\newcommand{\vota}{ V\ot A}
\newcommand{\vi}{V_i}
\newcommand{\vj}{V_j}
\newcommand{\dimcd}{\dim(\cd)}
\newcommand{\alij}{{\al_{ _{ij}}}}
\newcommand{\alji}{{\al_{ _{ji}}}}
\newcommand{\rcc}{r_\cc}
\newcommand{\rcd}{r_\cd}
\newcommand{\clsx}{[X]}
\newcommand{\clsy}{[Y]}
\newcommand{\clsz}{[Z]}
\newcommand{\rcdp}{r_{\cd'}}
\newcommand{\sumjtorp}{\sum_{j=0}^{r'}}
\newcommand{\aljm}{{\al_{ _{jm}}}}
\newcommand{\aljn}{{\al_{ _{jn}}}}
\newcommand{\sjm}{s_{jm}}
\newcommand{\smj}{s_{mj}}
\newcommand{\snj}{s_{nj}}
\newcommand{\betaij}{\beta_{ _{ij}}}
\newcommand{\betaji}{\beta_{ _{ji}}}
\newcommand{\gammaij}{\gamma_{ _{ij}}}
\newcommand{\gammaji}{\gamma_{ _{ji}}}
\newcommand{\ip}{i'}
\newcommand{\sumjtoprp}{\sum_{j=0}^{r'}}
\newcommand{\sumjtopr}{\sum_{j=0}^{r}}
\newcommand{\teh}{\tilde{h}}
\newcommand{\cdp}{{\cd'}}\newcommand{\xphii}{X_{\phi(i)}}
\newcommand{\inv}{^{-1}}
\newcommand{\fq}{\mtr f_{ Q}}
\newcommand{\tr}{\mtr{tr}}
\newcommand{\rtwone}{R_{21}R}
\newcommand{\ccad}{{\cc_{\mtr{ad}}}}
\newcommand{\ccpt}{{\cc_{\mtr{pt}}}}
\newcommand{\qtr}{quasi-triangular\;}
\newcommand{\trq}{\tr_q}
\newcommand{\repal}{\mtr{Rep}(A//L)}
\newcommand{\lkeravi}{\lker_A(V_i)}
\newcommand{\lkeravj}{\lker_A(V_j)}
\newcommand{\cross}[1][1pt]{\ooalign{%
\rule[1ex]{1ex}{#1}\cr
\hss\rule{#1}{.7em}\hss\cr}}
\newcommand{\blml}{\blam_L} 
\newcommand{\phir}{\phi_R}
\newcommand{\kda}{{  \Phi(A)}}
\newcommand{\mtcil}{\mtc{I}_L}
\newcommand{\un}{\unu}
\newcommand{\tfl}{\mtc{T}}
\newcommand{\barzm}{\barz(M)}
\newcommand{\barzn}{\barz(N)}
\newcommand{\ccr}{\mtc R^{\cc}}
\newcommand{\ulc}{\ul{\cc}}
\newcommand{\pimx}{\pi_{M;\;X}}
\newcommand{\pinx}{\pi_{N;\;X}}
\newcommand{\acc}{{\mathrm A_\cc}}
\newcommand{\epsu}{\eps_\unu}
\newcommand{\ob}{\mtr{Obj}}
\newcommand{\obc}{\mtr{Obj(\cc)}}
\newcommand{\ccop}{\cc^{\mtr{op}}}
\newcommand{\mtf}{\mtc F_\lam}
\newcommand{\mtfi}{\mtc F^{-1}_\lam}
\newcommand{\elcd}{\ell_\cd}
\newcommand{\mcid}{\mtc I_\cd}
\newcommand{\mcidp}{\mtc I_{\cd'}}
\newcommand{\wtildelcd}{\widetilde{\elcd}}
\newcommand{\wtildelcdp}{\widetilde{\ell_{\cd'}}}
\newcommand{\cpt}{\cc_{\mtr{pt}}}
\newcommand{\barzr}{\barz_\cd}
\newcommand{\barzv}{\barz(V)}
\newcommand{\acd}{\mathrm A_\cd}
\newcommand{\czrcd}{\cz_\cc(\cd)}
\newcommand{\sml}{\Small}
\newcommand{\bs}{\blue{\Small }}
\newcommand{\yd}{Yetter-Drinfeld}
\newcommand{\sumitor}{\sum_{i=0}^r}
\newcommand{\cdop}{\cd^{\mtr{op}}}
\newcommand{\ccrev}{\cc^{\mtr{rev}}}
\newcommand{\barz}{{\bar{\mathrm Z}}}
\newcommand{\etl}{etale\;}
\newcommand{\czca}{\cz(\ca)}
\newcommand{\tetx}{\text}
\newcommand{\widehta}{\widehat}
\newcommand{\wdhat}{\widehat}
\newcommand{\wht}{\widehat}
\newcommand{\cofa}{{\mathbb C[\mtc B]}}
\newcommand{\wdt}{\widehat}
\newcommand{\dl}{{^\#}}
\newcommand{\comx}{\mathbb C}
\newcommand{\sgj}{{\sg(j)}}
\newcommand{\mujo}{\mu_\jo}
\newcommand{\mujtw}{\mu_\jtw}
\newcommand{\adz}{a^{\#}}
\newcommand{\bdz}{b^{\#}}
\newcommand{\spr}{S^\perp}
\newcommand{\cofs}{\comp [S]}
\newcommand{\spz}{S^{\perp_z}}
\newcommand{\omz}{\omega_z}
\newcommand{\zg}{\mathrm{Z}(S)}
\newcommand{\aling}{{\al \in g}}
\newcommand{\blkg}{\mtr{Bl}(g)}
\newcommand{\clsg}{\mtr{Cl}(g)}
\newcommand{\mtadinv}{\mtc G^{{-1}}}
\newcommand{\muk}{\mu_{k}}
\newcommand{\mta}{\mtc F}
\newcommand{\cofad}{\comp[\wdht A]}
\newcommand{\wtau}{\wdht{\tau}}
\newcommand{\mtainv}{{\mta}^{-1}}
\newcommand{\wdht}{\widehat}
\newcommand{\augm}{\mtr{aug}}
\newcommand{\mua}{\wdht {\wdht a}}
\newcommand{\aps}{A//S}
\newcommand{\cssa}{\cc(S, A)}
\newcommand{\aug}{\mtr{aug}}
\newcommand{\rss}{{\big|_S}}
\newcommand{\gprp}{g^\perp}
\newcommand{\alins}{{s \in S}}
\newcommand{\sz}{s^{D}}
\newcommand{\wmu}{\widehta{\mu}}
\newcommand{\wmui}{\widehta{\mu}_i}
\newcommand{\wmuj}{\widehta{\mu}_j}
\newcommand{\wch}{\widehta{\ch}}
\newcommand{\wzd}{\widehat{d}}
\newcommand{\wpm}{\widehat{P}}
\newcommand{\wps}{\widehat{p}}
\newcommand{\gal}{\mtr{Gal}}
\newcommand{\galkq}{\gal(\mathbb K/\mathbb Q)}
\newcommand{\sgh}{\sg_{ _{H}}}
\newcommand{\sggi}{{\sg(i)}}
\newcommand{\sge}{\sg_{_{\widehat R}}}
\newcommand{\unue}{{\unu_{\cecc}}}
\newcommand{\mtcf}{\mtc {F}}
\newcommand{\wsgf}{\widehat{{\sg}_{ _F}}}
\newcommand{\sghstar}{{{\sg}_{ _{H^*}}}}
\newcommand{\we}{\widehta{E}}
\newcommand{\sumktom}{\sum_{k=0}^m}
\newcommand{\wf}{\widehat{F}}
\newcommand{\hsgj}{\widehat{\sg}(j)}
\newcommand{\whsgi}{\widehta{\sg}(i)}
\newcommand{\wpp}{\widehat{p}}
\newcommand{\tauj}{{\tau(j)}}
\newcommand{\dimcctauj}{\dim(\cc^\tauj)}
\newcommand{\etas}{{\eta(s)}}
\newcommand{\mcc}{m_H}
\newcommand{\wal}{\widehta{\al}}
\newcommand{\wj}{\widehat{\mtc J}}
\newcommand{\galc}{\mtr{Gal}_{\cc}}
\newcommand{\galz}{\mtr{Gal}_{\czcc}}
\newcommand{\wjr}{\widehat{J}_{R}}
\newcommand{\dimcck}{\dim(\cc^k)}
\newcommand{\wgrcc}{\widehat{\mtr{Gr}(\cc)}}
\newcommand{\nchi}{{\frac{\ch_i}{\di }}} 
\newcommand{\nchj}{{\frac{\ch_j}{\dxj}}}
\newcommand{\wni}{{\widehat{n}_i}}
\newcommand{\sgte}{\widetilde{\sg_E}}
\newcommand{\mtad}{\mtc G}
\newcommand{\whj}{\widehta{h}_j}
\newcommand{\jdl}{{j\dl}}
\newcommand{\wcfcc}{\widehat{\cfcc}}
\newcommand{\mutauj}{\mu_{\tau(j)}}
\newcommand{\tauk}{\tau(k)}
\newcommand{\muzm}{{\mu_0^{-}}}
\newcommand{\sqrtog}{\sqrt{|G|}}
\newcommand{\muz}{\mu_0}
\newcommand{\njtw}{n_\jtw}
\newcommand{\njo}{n_\jo}
\newcommand{\fjo}{F_\jo}
\newcommand{\fjtw}{F_\jtw}
\newcommand{\wta}{\widehat{A}}
\newcommand{\dol}{{^{\circ}}}
\newcommand{\bdl}{{b\dl}}
\newcommand{\jdol}{{j\dol}}
\newcommand{\fj}{F_j}
\newcommand{\cwta}{\comp[\wta]}
\newcommand{\hx}{\widehta{x}}
\newcommand{\hy}{\widehta{y}}
\newcommand{\cal}{\mtc A_{\al}}
\newcommand{\xuu}{x_{uu}}
\newcommand{\wxuu}{\widehat{\xuu}}
\newcommand{\xvv}{x_{vv}}
\newcommand{\xuv}{x_{uv}}
\newcommand{\xmn}{x_{m,n}}
\newcommand{\buvmn}{B^{u,v}_{m,n}}
\newcommand{\blm}{\blam}
\newcommand{\dimccr}{\dim(\cc^r)}
\newcommand{\adl}{a\dl}
\newcommand{\sumltom}{\sum_{l=0}^m}
\newcommand{\mbq}{\mathbb Q}
\newcommand{\mbqs}{\mathbb Q(S)}
\newcommand{\mbk}{\mathbb K}
\newcommand{\mz}{\mathbb Z}
\newcommand{\wsgj}{\widehat{\sigma}(j)}
\newcommand{\wsgi}{\widehat{\sigma}(i)}
\newcommand{\wg}{\widehat{g}}
\newcommand{\wtf}{\widehat{F}}
\newcommand{\galqspq}{\mtr{Gal}(\mathbb Q(S)/\mathbb Q)}
\newcommand{\cctauj}{\cc^{\tau(j)}}
\newcommand{\cctauk}{\cc^{\tau(k)}}
\newcommand{\wtfj}{\widetilde{F_j}}
\newcommand{\wfj}{\widetilde{F_j}}
\newcommand{\wtmuj}{\widetilde{\mu_j}}
\newcommand{\wmtcfj}{{\widetilde{\mtc F}_j}}
\newcommand{\mtfr}{\mtr{F_a}}
\newcommand{\wdr}{R_\comp^*}
\newcommand{\fgph}{{F_{G/H}}}
\newcommand{\wcfj}{\wmtcfj}
\newcommand{\nxi}{{\frac{x_i}{\di }}}
\newcommand{\fpr}{{\fp(R)}}
\newcommand{\nxs}{{\frac{x_s}{d_s}}}
\newcommand{\mtfme}{\mtc F}
\newcommand{\chic}{\ch_i^{\circ}}
\newcommand{\chjc}{\ch_j^{\circ}}
\newcommand{\mtfsh}{{\mtc F_\lam}}
\newcommand{\mupq}{{\mu_{pq}}}
\newcommand{\tlam}{{\widetilde{\lam}}}
\newcommand{\chid}{{\ch_i^{\circ}}}
\newcommand{\rc}{{R_\comp}}
\newcommand{\rgo}{{\mathbb R_{\geq 0}}}
\newcommand{\sumrorc}{{{\sum\limits_{\ro \in \rc}}}}
\newcommand{\aliro}{{\ro(x_i)}}
\newcommand{\barjd}{{\bar{\mtc J_\cd}}}
\newcommand{\lbarcj}{{\frac{C_j}{{\dim(\mathcal C^j)}}}}
\newcommand{\omtcb}{{\overline{\mathcal B}}}
\newcommand{\whr}{{\widehat{R}}}
\newcommand{\nxj}{\frac{x_j}{\dxj}}
\newcommand{\nxk}{\frac{x_k}{d_k}}
\newcommand{\onkij}{{\overline{N^k_{ij}}}}
\newcommand{\sgk}{{\sigma(k)}}
\newcommand{\sgl}{{\sigma(l)}}
\newcommand{\fqi}{{\fq^{-1}}}
\newcommand{\wdb}{{\widehat{\mtc B}}}
\newcommand{\mtcb}{{\mtc B}}
\newcommand{\nif}{{h_i}}
\newcommand{\rb}{(R, \mtc B)}
\newcommand{\nxip}{{\frac{x_{i'}}{d_{i'}}}}
\newcommand{\mujp}{{\mu_{j'}}}
\newcommand{\etai}{{\eta(i)}}
\newcommand{\wsg}{{\widehat{\sg}}}
\newcommand{\wsgh}{{\wsg_{ _{H}}}}
\newcommand{\wtaui}{{\widehat{\tau}(i)}}
\newcommand{\wsghstar}{{\wsg_{H^*}}}
\newcommand{\wtauj}{{\wtau(j)}}
\newcommand{\weta}{{\widehat{\eta}}}
\newcommand{\detai}{{d_{ _{\eta(i)}}}}
\newcommand{\hbz}{{(H, \mtc B, \mu_0)}}
\newcommand{\whbz}{{\widehat{\hbz}}}
\newcommand{\tsgh}{{{\widetilde{\sgh}}}}
\newcommand{\ghb}{{G_\hb}}
\newcommand{\taujo}{{\tau(\jo)}}
\newcommand{\taujtw}{{\tau(\jtw)}}
\newcommand{\mutauk}{{\mu_{\tauk}}}
\newcommand{\hetai}{{h_{ _\etai}}}
\newcommand{\xetai}{{x_{ _\etai}}}
\newcommand{\wn}{{\widehat{n}}}
\newcommand{\wh}{{\widehat{h}}}
\newcommand{\distar}{{d_{i^*}}}
\newcommand{\dwtaui}{{d_{ _{\wtaui}}}}
\newcommand{\sumiptom}{{\sum_{\ip=0}^m}}
\newcommand{\alitaugj}{{\al_{ _{i\tau_g(j)}}}}
\newcommand{\dtaugj}{{d_{ _{\tau_g(j)}}}}
\newcommand{\mip}{{M(\ip)}}
\newcommand{\sumttom}{{\sum_{t=0}^m}}
\newcommand{\muxi}{{\mu_{ _{[X_i]}}}}
\newcommand{\muxip}{{\mu_{ _{[X_{\ip}]}}}}
\newcommand{\ncj}{{\frac{C_j}{\dim(\cc^j)}}}
\newcommand{\jp}{{j'}}
\newcommand{\minv}{{M^{-1}}}
\newcommand{\tfq}{{\widehat{\fq}}}
\newcommand{\catcecc}{{\mtr{CE}(\cc)}}
\newcommand{\wcatfcc}{{\widehat{\catfcc}}}
\newcommand{\mforall}{{\;\;\text{for all}\;\;}}
\newcommand{\what}{\widehat}
\newcommand{\wfz}{{\widehat{F}_0}}
\newcommand{\hbfr}{{(H, \mtc B, \fp)}}
\newcommand{\sgn}{{\mtr{sgn}}}
\newcommand{\wir}{{\widehat{R}}}
\newcommand{\gcc}{{G(\cc)}}
\newcommand{\jccpt}{{J_{ _{\ccpt}}}}
\newcommand{\jccad}{{J_{ _{\ccad}}}}
\newcommand{\fpccad}{{\fp(\ccad)}}
\newcommand{\fpccpt}{{\fp(\ccpt)}}
\newcommand{\kc}{{K(\cc)}}
\newcommand{\wkc}{{\widehat{\kc}}}
\newcommand{\hs}{{(L,\;\mtc S)}}
\newcommand{\rbad}{{H_{ _{ad}}}}
\newcommand{\hbad}{{\hb_{ad}}}
\newcommand{\jhbad}{{J_{\hbad}}}
\newcommand{\htt}{{(K, \mtc T)}}
\newcommand{\jhtt}{{\mtc J_{ _{\htt}}}}
\newcommand{\jhs}{{\mtc J_{ _{\hs}}}}
\newcommand{\lamhs}{{\lam_{ _{\hs}}}}
\newcommand{\lamhtt}{{\lam_{ _{\htt}}}}
\newcommand{\coo}{{co}}
\newcommand{\wdhad}{{(\wdh)_{ad}}}
\newcommand{\had}{{H_{ad}}}
\newcommand{\wdh}{\widehat{H}}
\newcommand{\whbad}{{\whb_{ _{ad}}}}
\newcommand{\kerhb}{{\ker_{ _{\hb}}}}
\newcommand{\gwdh}{{G(\wdh)}}
\newcommand{\nxl}{\frac{x_l}{d_l}}
\newcommand{\proditom}{{\prod_{i=0}^m}}
\newcommand{\wdhn}{{\wdh^{(n)}}}
\newcommand{\hn}{{H_{(n)}}}
\newcommand{\nox}{{\fp(x)}}
\newcommand{\mujstar}{{\mu_{j^\#}}}
\newcommand{\sco}{{S^\coo}}
\newcommand{\rrad}{{I(1)}}
\newcommand{\istar}{{i^*}}
\newcommand{\mtcs}{{\mtc S}}
\newcommand{\lamccad}{{{\lam_{\ccad}}}}
\newcommand{\lamczccad}{{{\lam_{\cczc_{ad}}}}}
\newcommand{\wtm}{{\widetilde{m}}}
\newcommand{\cczerog}{{\cc^0_G}}
\newcommand{\zd}{{\mathbb Z_d}}
\newcommand{\ucc}{{U(\cc)}}
\newcommand{\uczcc}{{U(\czcc)}}
\newcommand{\czccpt}{{\czcc_{pt}}}
\newcommand{\cg}{{[G,G]}}
\newcommand{\fjss}{F^j_{ss}}
\newcommand{\csujss}{\mtr{C}^j_{ss}}
\newcommand{\csujts}{\mtr{C}^j_{ts}}
\newcommand{\fjts}{F^j_{ts}}
\newcommand{\iotal}{\iota_{\mthl}}
\newcommand{\hmcd}{{\hm_\cd}}
\newcommand{\csuj}{\csu^j}
\newcommand{\mujst}{\mu^j_{st}}
\newcommand{\mukuv}{\mu^k_{uv}}
\newcommand{\mujpsptp}{\mu^{j'}_{s't'}}
\newcommand{\cjst}{\csu^j_{st}}
\newcommand{\cjts}{\csu^j_{ts}}
\newcommand{\cjpsptp}{\csu^{j'}_{s't'}}
\newcommand{\cjptpsp}{\csu^{j'}_{t's'}}
\newcommand{\ncjst}{\frac{\cjst}{\dim(\cc_j)}}
\newcommand{\ncjts}{\frac{\cjts}{\dim(\cc_j)}}
\newcommand{\ncjpsptp}{\frac{\cjpsptp}{\dim(\cc_{j'})}}
\newcommand{\ncjptpsp}{\frac{\cjptpsp}{\dim(\cc_{j'})}}
\newcommand{\lu}{\mtc L(\unu)}
\newcommand{\lmod}{L-\mtr{mod}}
\newcommand{\rex}{\mtr{Rex}}
\newcommand{\resml}{\mtr{res}^A_L}
\newcommand{\cbt}{\;\underline{{\cdot}}\;}
\newcommand{\acm}{{A_{\cm}}}
\newcommand{\cfm}{\cf(\cm)}
\newcommand{\res}{{\mtr{Res}}}
\newcommand{\repk}{{\rep(K)}}
\newcommand{\resra}{{\res^{\mtr{ra}}}}
\newcommand{\rha}{{\mtr{\ra}}}
\newcommand{\epsul}{\eps_L}
\newcommand{\epsm}{\eps_M}
\newcommand{\acml}{(\al^c_M)|_L}
\newcommand{\idm}{\id_M}
\newcommand{\idl}{\id_L}
\newcommand{\deltau}{\delta_{\unu}}
\newcommand{\ptr}{\mtr{ptr}}
\newcommand{\zom}{{Z(M)}}
\newcommand{\csl}{{\cs_L}}
\newcommand{\ztmn}{{Z^{(2)}(M,N)}}
\newcommand{\bfastar}{\bfa^*}
\newcommand{\zu}{{Z(\unu)}}
\newcommand{\di}{d_i}
\newcommand{\ccro}{{\cc^{\rho}}}
\newcommand{\ccpj}{{\cc^{(j)}}}
\newcommand{\cfa}{\rc}
\newcommand{\cflo}{\cf(L_1)}
\newcommand{\cfltw}{\cf(L_2)}
\newcommand{\fjst}{F^{j}_{st}}
\newcommand{\ccpsj}{{\cc^{(j)}_s}}
\newcommand{\ccpsro}{{\cc^{\ro}_s}}
\newcommand{\ccpjs}{\ccpsj}
\newcommand{\ccjps}{\ccpsj}
\newcommand{\ccjsp}{\ccpsj}
\newcommand{\tfjst}{\widetilde{\fjst}}
\newcommand{\lj}{\mtc L_j}
\newcommand{\mtmj}{\mtc M_j}
\newcommand{\tfjpsptp}{{\widetilde{F^{j'}_{s't'}}}}
\newcommand{\tfjstp}{{\widetilde{F^{j}_{st'}}}}
\newcommand{\enxczcca}{{\enx_\czcc(A)}}
\newcommand{\ccjpt}{\cc^{(j)}_t}
\newcommand{\ccpjt}{\cc^{(j)}_t}
\newcommand{\sumjtomst}{\sum_{j\in \mtc J}\;\sum_{s,t\in \mtmj}}
\newcommand{\fjpsptp}{{F^{j'}_{s't'}}}
\newcommand{\fjstp}{{F^{j}_{st'}}}
\newcommand{\alijst}{{\al(i)^j_{st}}}
\newcommand{\csujst}{\csu^j_{st}}
\newcommand{\st}{_{st}}
\newcommand{\mtha}{\mathbb A}
\newcommand{\bcl}{\bar{\mtr{CE}}(L)}
\newcommand{\jcd}{\mtc J_{\cd}}
\newcommand{\jce}{\mtc J_{\ce}}
\newcommand{\ls}{{(L, S)}}
\newcommand{\jtwo}{j_{2}}
\newcommand{\hmda}{\hm_{D(A)}}
\newcommand{\uwha}{\uw_\lsh^A}
\newcommand{\mci}{{\mtc I}}
\newcommand{\irrda}{\irr(D(A))}
\newcommand{\lsg}{{L(g)}}
\newcommand{\lsh}{{L(h)}}
\newcommand{\dimkm}{\dimk(M)}
\newcommand{\lsm}{{L(m)}}
\newcommand{\phia}{\phi(A)}
\newcommand{\jcl}{\mtc J_L}
\newcommand{\kastar}{{K(A^*)}}
\newcommand{\cajojtw}{\ca(j_1, j_2)}
\newcommand{\mtcj}{\mtc J}
\newcommand{\cbjojtw}{\ccb(j_1, j_2)}
\newcommand{\nioitw}{N^i_{\io, \itw}}
\newcommand{\barcmio}{{\bar C}_{m(\io)}}
\newcommand{\barcmitw}{{\bar C}_{m(\itw)}}
\newcommand{\alio}{\al_{ _\io}}
\newcommand{\ali}{\al_{ _i}}
\newcommand{\alitw}{\al_{ _\itw}}
\newcommand{\barci}{\bar C_{i}}
\newcommand{\barcmi}{\bar C_{m(i)}}
\newcommand{\mio}{m(\io)}
\newcommand{\mitw}{m(\itw)}
\newcommand{\ccmio}{\cc^{m(\io)}}
\newcommand{\ccmi}{\cc^{m(i)}}
\newcommand{\ccmitw}{\cc^{m(\itw)}}
\newcommand{\ccjo}{\cc^{\jo}}
\newcommand{\ccjtw}{\cc^{\jtw}}
\newcommand{\cscm}{\cs_\cm}
 \newcommand{\wloge}{without loss of generality }
\newcommand{\barzmu}{ \barz_{\cm}(\unu)}
\newcommand{\reg}{\mtr{Reg}}
\newcommand{\barzrm}{ \barz_{\cm}}
\newcommand{\esuv}{F^s_{uv}}
\newcommand{\esvu}{F^s_{vu}}
\newcommand{\esuu}{F^s_{uu}}
\newcommand{\esvv}{F^s_{vv}}
\newcommand{\tildef}{\tilde F}
\newcommand{\ovl}{\overline}
\newcommand{\sccj}{{|\cc^j|}}
\newcommand{\dimccj}{\dim(\cc^j)}
\newcommand{\ccl}{\cc^l}
\newcommand{\cck}{\cc^k}
\newcommand{\ccp}{{\cc'}}
\newcommand{\gr}{\mtr{Gr}}
\newcommand{\mtlsum}{\mathlarger{\sum}}
\newcommand{\whag}{{\widehat{(a, \gamma)}}}
\newcommand{\grcc}{\gr(\cc)}
\newcommand{\sumstom}{\sum_{s=0}^m}
\newcommand{\barcr}{\bar C_r}
\newcommand{\barcl}{\bar C_l}
\newcommand{\barcj}{\bar C_j}
\newcommand{\barck}{\bar C_k}
\newcommand{\io}{{i_1}}
\newcommand{\itw}{{i_2}}
\newcommand{\jo}{{j_1}}
\newcommand{\jtw}{{j_2}}
\newcommand{\barcjo}{\bar C_{j_1}}
\newcommand{\barcjtw}{\bar C_{j_2}} 
\newcommand{\chitw}{\chi_{\itw}}
\newcommand{\xitw}{x_{ _{\itw}}}
\newcommand{\mi}{M(i)}
\newcommand{\chio}{\ch_{\io}}
\newcommand{\xio}{x_{ _{\io}}}
\newcommand{\phibr}{\phi_{\bar R}}
\newcommand{\jc}{{\mtc J}^{c}}
\newcommand{\cvi}{C_{V_i}}
\newcommand{\tausg}{{\tau_{ _\sg}}}
\newcommand{\cdl}{\cd_L}
\newcommand{\cfl}{\overline{\mtr{CF}}(L)}
\newcommand{\cel}{\overline{\mtr{CE}}(L)}
\newcommand{\mtcfj}{\mtc F_j}
\newcommand{\jl}{\mtc J_L}
\newcommand{\resal}{\mtr{Res}_L}
\newcommand{\lamcdl}{\lam_{\cdl}}
\newcommand{\mmj}{\mtc M_L}
\newcommand{\zim}{\{0,\dots , m\}}
\newcommand{\zir}{\{0,\dots , r\}}
\newcommand{\cdm}{\cd_M}
\newcommand{\cdn}{\cd_N}
\newcommand{\jm}{\mtc J_M}
\newcommand{\resalfj}{\resal(F_j)}
\newcommand{\zh}{\mtc Z(H)}
\newcommand{\enxcza}{\enx_{\czcc}(A)}
\newcommand{\bcj}{\bar C_j}
\newcommand{\bpcj}{\bar C'_j}
\newcommand{\lamcdm}{\lam_{\cd_M}}
\newcommand{\lamcdn}{\lam_{\cd_N}}
\newcommand{\ena}{\enx_\czcc(A)}
\newcommand{\sumutom}{{\sum_{u=0}^m}}
\newcommand{\mti}{\mtc I}
\newcommand{\svec}{\mtr{sVec}}
\newcommand{\irrcc}{\irr(\cc)}
\newcommand{\irrcd}{\irr(\cd)}
\newcommand{\cocd}{\co(\cd)}
\newcommand{\cocc}{\co(\cc)}
\newcommand{\sgst}{\sg^{\mtr{st}}}
\newcommand{\sgstzv}{\sg^{\mtr{st}}_{Z(V), -}}
\newcommand{\dltv}{\delta_V}
\newcommand{\bdlt}{{\bar{\delta}}}
\newcommand{\bpivx}{\bar{\pi}_{V;\;X}}
\newcommand{\bpimx}{\bar{\pi}_{M;\;X}}
\newcommand{\bpinx}{\bar{\pi}_{N;\;X}}
\newcommand{\bpi}{\bar{\pi}}
\newcommand{\mtj}{\mtc L}
\newcommand{\mtjcd}{\mtj_{\cd}}
\newcommand{\lag}{\langle}
\newcommand{\rag}{\rangle}
\newcommand{\dvx}{\delta_{V, X}}
\newcommand{\evx}{{\ev_X}}
\newcommand{\barzw}{{\barz(W)}}
\newcommand{\bphi}{\bar{\phi}}
\newcommand{\zcat}{Z_{\mtr{cat}}}
\newcommand{\repdh}{\rep(D(H)}
\newcommand{\rrange}{R()}
\newcommand{\lrange}{L()}
\newcommand{\dimkh}{\dimk(H)}
\newcommand{\reph}{\rep_{\kk}(H)}
\newcommand{\omg}{\Omega}
\newcommand{\hdl}{\hat{\mtc D}_L}
\newcommand{\phiro}{{\phi_\ro}}
\newcommand{\tev}{\tilde{ev}}
\newcommand{\betai}{\beta^{-1}}
\newcommand{\alphai}{\alpha^{-1}}
\newcommand{\uu}{{(1,1)}}
\newcommand{\tu}{{(2,1)}}
\newcommand{\td}{{(2,2)}}
\newcommand{\ud}{{(1,2)}}
\newcommand{\dtwo}{{(2)}}
\newcommand{\bara}{{\bar{A}}}
\newcommand{\minu}{{(-1)}}
\newcommand{\zero}{{(0)}}
\newcommand{\fic}{F_i^c}
\newcommand{\chc}{\ch^c}
\newcommand{\barga}{\bar G(A)}
\newcommand{\repa}{\rep(A)}
\newcommand{\repaad}{\rep(A)_{\mtr{ad}}}
\newcommand{\repbad}{\rep(B)_{\mtr{ad}}}
\newcommand{\repb}{\rep(B)}
\newcommand{\rp}{\rep}
\newcommand{\cmhhh}{_H\cm_H^H}
\newcommand{\cmhh}{_H\cm_H}
\newcommand{\hot}{\hat{\ot}}
\newcommand{\bars}{\bar S}
\newcommand{\cmh}{_H\cm}
\newcommand{\phiu}{\phi_1}
\newcommand{\bphiu}{\bar{\phi}_1}
\newcommand{\phit}{\phi_2}
\newcommand{\bphit}{\bar{\phi}_2}
\newcommand{\phitr}{\phi_3}
\newcommand{\bphitr}{\bar{\phi}_3}
\newcommand{\qur}{q_1^{\mtr{R}}}
\newcommand{\qtwor}{q_2^{\mtr{R}}}
\newcommand{\coih}{{\mtr{co} \; H}}
\newcommand{\qhb}{\text{quasi-Hopf bimodules}\;}
\newcommand{\wtp}{{\widehat{P}}}
\newcommand{\qtwol}{q_2^{\mtr L}}
\newcommand{\qul}{q_1^{\mtr L}}
\newcommand{\hcmhh}{\;^H_H\cm_H}
\newcommand{\rphi}{{\;_R\phi}}
\newcommand{\tblam}{{\tilde{\blam}}}
\newcommand{\qu}{{q_\unu}}
\newcommand{\fpch}{\fp(\ch)}
\newcommand{\sgi}{\sigma(i)}
\newcommand{\ccrbz}{\ccr^{\barz}}
\newcommand{\mtl}{{  \mtc L}}
\newcommand{\mtlu}{{  \mtc L}(\unu)}
\newcommand{\piuy}{\pi_{\unu, Y}}
\newcommand{\piux}{\pi_{\unu, X}}
\newcommand{\iuy}{\iota_{\unu, Y}}
\newcommand{\iux}{\iota_{\unu, X}}
\newcommand{\rhom}{\rho_M}
\newcommand{\z}{{  Z}}
\newcommand{\czx}{\cz(X)}
\newcommand{\czm}{\z(M)}
\newcommand{\czv}{\cz(V)}
\newcommand{\zm}{\cz(M)}
\newcommand{\pivx}{\pi_{V;\;X}}
\newcommand{\barza}{\barz(A)}
\newcommand{\hra}{\hookrightarrow}
\newcommand{\barj}{\mathbf J}
\newcommand{\mul}{\mu_{\barzm}^l}
\newcommand{\mur}{\mu_{\barzm}^r}
\newcommand{\chofm}{\ch_M}
\newcommand{\pium}{\pi_{\unu, M}}
\newcommand{\prc}{\perp^{\mtr{proj}}}
\newcommand{\wrt}{with respect to\;}
\newcommand{\fpas}{\frac{p_s}{a_s} }
\newcommand{\sgs}{{\sigma(s)}}
\newcommand{\ccsgs}{\cc_{\sgs}}
\newcommand{\cu}{\mtc U}
\newcommand{\czzcc}{\cz(\czcc)}
\newcommand{\cci}{\cc^i}
\newcommand{\cfzcc}{{\cf(\czcc)}}
\newcommand{\czccad}{\czcc_{\ad}}
\newcommand{\vsk}{\vskip 0.15cm \noindent}
\newcommand{\ccj}{\cc^j}
\newcommand{\chad}{\ch_{\ad}}
\newcommand{\sumstor}{\sum_{s=0}^r}
\newcommand{\fpcc}{\fp(\cc)}
\newcommand{\fpx}{\fp(X)}
\newcommand{\ccs}{\cc^s}
\newcommand{\cct}{\cc^t}
\newcommand{\ircc}{\irr(\cc)}
\newcommand{\ccsuu}{\ccs_{uu}}
\newcommand{\lamad}{\lam_{\ad}}
\newcommand{\regcc}{R_\cc}
\newcommand{\gch}{G_{\ch}}
\newcommand{\gmu}{G_\mu}
\newcommand{\mus}{\mu_s}
\newcommand{\gze}{G_0}
\newcommand{\chun}{\ch(1)}
\newcommand{\pkjojtw}{{\wpp_k(\jo, \jtw)}}
\newcommand{\pccj}{p_{\cc^j}}
\newcommand{\wtpccj}{\widetilde{\pccj}}
\newcommand{\bxd}{\boxed}
\newcommand{\incl}{\hookrightarrow}
\newcommand{\lha}{^A\lh}
\newcommand{\fbx}{\fbox}
\newcommand{\tsig}{\tilde\sigma}
\newcommand{\mir}{{\mtr{mir}}}
\newcommand{\wtl}{\widetilde}
\newcommand{\forg}{\mtr{Forg}}
\newcommand{\barh}{\underline H}
\newcommand{\alx}{\al_X}
\newcommand{\bxt}{\boxtimes}
\newcommand{\barhx}{\barh(X, X)}
\newcommand{\barcx}{\barc(X, X)}
\newcommand{\czcm}{\cz(\cc^*_\cm)}
\newcommand{\cstrm}{\cc^*_\cm}
\newcommand{\apk}{A//K}
\newcommand{\stc}{^{*\cop}}
\newcommand{\ka}{\mtr{K}(A)}
\newcommand{\xij}{x_{ij}}
\newcommand{\xii}{x_{ii}}
\newcommand{\barg}{{  \overline{G}}}
\newcommand{\vect}{\mtr{Vec}}
\newcommand{\jcdu}{J_{\cd}(\unu)}
\newcommand{\sent}{\mapsto}
\newcommand{\xj}{X_j}
\newcommand{\mui}{\mu_i}\newcommand{\muj}{\mu_j}
\newcommand{\kercc}{{\ker_{\cc}}}
\newcommand{\kerzcc}{{\ker_{\czcc}}}
\newcommand{\lkercc}{\lker_\cc}
\newcommand{\lkercci}{\lker_\cc(X_i)}
\newcommand{\tildefj}{\tilde{F}_j}
\newcommand{\fstar}{F_*}
\newcommand{\bargcc}{\barg(\cc)}
\newcommand{\tildecd}{\tilde{\cd}}
\newcommand{\kcc}{\mtr{K}(\cc)}
\newcommand{\sgpi}{\sg_\pi}
\newcommand{\barzx}{\barz(X)}
\newcommand{\ccrm}{\ccr_M}
\newcommand{\ccrn}{\ccr_M}
\newcommand{\barc}{{{ \ul{ \cc}}}}
\newcommand{\hmc}{\hm_\cc}
\newcommand{\rdx}{^*X}
\newcommand{\excl}{\blue{\bf !!!!!!!}}
\newcommand{\gs}{\mapsto}
\newcommand{\gos}{\mapsto}
\newcommand{\barzmn}{\barz(M\ot N)}
\newcommand{\hkw}{\hookrightarrow}
\newcommand{\qst}{q_*}
\newcommand{\qust}{q^*}\newcommand{\cec}{\mtr{  CE}(\cc)}
\newcommand{\ced}{\mtr{  CE}(\cd)}
\newcommand{\barzd}{\barz_{\cd\hookrightarrow \cc}}
\newcommand{\bll}{\blue}
\newcommand\QQ{{\mathbb Q}} 
\newcommand\RR{{\mathbb R}} 
\newcommand\CC{{\mathbb C}} 
\newcommand\FF{{\mathbb F}}
\renewcommand\SS{\mathfrak{S}}
\newcommand\htimes{\,\widehat\otimes\,}
\newcommand\pib{\overline{\pi}}  
\renewcommand\top{\mathrm{top}\,\!}
\newcommand\bS{\mathbf{S}}
\newcommand\C{\mathcal{C}}
\newcommand\groth{{G_0(A)}}
\newcommand\xx{\mathbf{x}}
\newcommand\sss{\mathbf{s}}
\newcommand\ppp{\mathbf{p}}
\newcommand\nnn{\mathbf{n}}
\newcommand\GG{\mathbf{G}}
\newcommand\NN{\mathbb{N}}
\newcommand{\ccb}{{\mathcal B}}
\newcommand{\bdelta}{\Delta}
\newcommand{\bann}{\mtr{Ann}}

\title[Conjugacy class dimensions]
{Frobenius-Perron Dimensions of Conjugacy Classes and an Ito-Michler-Type Result in Modular Fusion Categories}
\author{Sebastian Burciu}
\address{Inst.\ of Math.\ ``Simion Stoilow" of the Romanian Academy P.O. Box 1-764, RO-014700, Bucharest, Romania}
\email{sebastian.burciu@imar.ro}


\date{\today}
\subjclass{16T30; 18M20}
\bd
\thanks{The author is supported by a Research Grant GAR 2023 (cod 114), supported from the Donors' Recurrent Fund of the Romanian Academy,
managed by the "PATRIMONIU" Foundation.}

\maketitle
\begin{abstract}
The influence of certain arithmetic conditions on the sizes of conjugacy classes of a finite group on the group structure has been extensively studied in recent years. In this paper, we explore analogous properties for fusion categories. In particular, we establish an Ito-Michler-type result for modular fusion categories.
\end{abstract}
\setcounter{tocdepth}{1}
\newcommand{\wlamczccad}{{\widetilde{\lam}_{\czccad}}}
\newcommand{\wlamccad}{{\widetilde{\lam}_{\ccad}}}
\newcommand{\oucc}{{|\ucc|}}
\newcommand{\ccadpt}{{{\ccad}_{pt}}}
\newcommand{\oG}{{|G|}}
\newcommand{\oN}{{|N|}}
\newcommand{\oH}{{|H|}}
\newcommand{\ccgad}{{\cc^G_{ad}}}
\newcommand{\dimCi}{{\dim(\cc^i)}}
\newcommand{\dimCj}{{\dim(\cc^j)}}
\newcommand{\dimCk}{{\dim(\cc^k)}}
\newcommand{\ckij}{{c^k_{ij}}}
\newcommand{\pkij}{{p^k_{ij}}}
\newcommand{\Vi}{{V_i}}
\newcommand{\wtd}{\widetilde}
\newcommand{\mtrce}{{\mtr{CE}}}
\newcommand{\wtmu}{\widetilde{\mu}}
\newcommand{\ceczcc}{\mtrce(\czcc)}
\newcommand{\wtlam}{\widetilde{\lam}}
\newcommand{\fpcd}{{\fp(\cd)}}
\newcommand{\dimccpt}{{\dim(\ccpt)}}
\newcommand{\tch}{\widetilde{\ch}}
\newcommand{\chp}{{\ch'}}
\newcommand{\chpp}{{\ch^{''}}}
\newcommand{\mtcv}{{\mtc V}}
\newcommand{\xic}{{x_i^{\circ}}}
\newcommand{\xjc}{{x_j^{\circ}}}

\newcommand{\clcrd}{\mtr{Cl}(\cc:\cd)}
\newcommand{\sumjinc}{\sum_{j\in C}}
\newcommand{\pic}{\pi_e}
\newcommand{\piz}{\pi_e}
\newcommand{\pif}{\pi_f}
\newcommand{\pt}{{\mtr{pt}}}
\newcommand{\dab}{{\cdot_{\mtr{ad}}}}
\newcommand{\pij}{{\pi(j)}}
\newcommand{\prs}{p_{r,s}}
\newcommand{\lbardm}{\overline{d}_m}
\newcommand{\lbardt}{\frac{\mtr{D}_t}{\dim(\cd^t)}}
\newcommand{\pie}{{\pi_e}}

\newcommand{\psixi}{{\psi_{[X_i]}}}
\newcommand{\psixip}{{\psi_{[X_\ip]}}}
\newcommand{\lamcdp}{{\lam_\cdp}}
\section{Introduction}\label{intro}
Fusion categories have emerged as a powerful tool in both mathematics and theoretical physics, offering a categorical framework for understanding symmetry and structure in various contexts, such as quantum field theory, topology, and representation theory. These categories generalize many classical notions from group theory, especially in settings where group symmetries alone are insufficient to describe complex phenomena.

Over the past few decades, considerable attention has been given to understanding the structure of a finite group $G$ based on information about the arithmetic properties of its irreducible character degrees or, dually, of the sizes of its conjugacy classes. One of the central questions that has been extensively studied is how the structure of a group $G$ can be determined when certain arithmetic properties of its irreducible character set, $\text{Irr}(G)$, are known. Significant contributions to this area have been made recently.  For a summary of the known results in this direction one may consult the expository papers of \cite{Manz}  and \cite{lewis}.

 Building on the work of Zhu \cite{Z}, Cohen and Westreich \cite{CW6} introduced the concept of conjugacy classes for semisimple Hopf algebras. This new framework allowed them to extend several results from the theory of finite group representations to the context of semisimple Hopf algebra representations. More recently, Shimizu \cite{scalg} expanded on this idea by defining conjugacy classes for fusion categories, generalizing the notion introduced by Cohen and Westreich. In the same work, Shimizu associated a central element, termed the conjugacy class sum, to each conjugacy class, which parallels the role of the sum of elements in a group conjugacy class in group theory. Shimizu also posed the question of which results from \cite{CW6, CW4} could be further extended from semisimple Hopf algebras to fusion categories.

The primary goal of this article is to establish analogous results for the Frobenius-Perron dimensions of conjugacy classes in weakly integral fusion categories and to compare these findings with the corresponding results for finite groups. It is noteworthy that there are currently very few results concerning the Frobenius-Perron dimensions of conjugacy classes in fusion categories. Cohen and Westreich proved  in \cite[Theorem 4.5]{CW4} that a simple non-commutative quasitriangular semisimple Hopf algebra possesses no non-trivial conjugacy classes whose dimensions are powers of a prime number.

The key ingredient in our study is Harada's identity on the product of conjugacy class sums  for weakly integral fusion categories, see \cite[Theorem 1.2]{harada}.  Our first main results in this direction are the following:

Let $\cc$ be a weakly integral fusion category and  $\{\ccj\}_{j=0}^m$ the set of all its conjugacy classes and $C_j$ their associated class sums.
\bt\label{gen-dual}
In any weakly integral braided fusion category $\cc$  with a commutative Grothendieck ring $K(\cc)$ the following holds:
\beq\label{div:fpcc}
\frac{\fp(\ccpt)}{\fp(\cc)}\big[\prod_{j=0}^m\fp(\ccj)\big]^{2\rank(\cc)+1}\in \mathbb Z.
\eeq
\et
For the definition of conjugacy classes $\cc^j$ and their associated class sums $C_j$, refer to Section \ref{prelim}. Moreover, $\ccpt$ denotes the maximal pointed fusion subcategory of  $\cc$.
\bt \label{ccj-not div-p}
Let $\cc$ be a weakly integral braided fusion category, and let $p$ be a prime divisor of $\fpcc$ such that $\fpcc=p^\al N$ with $\al\geq 1$ and $(p, N)=1$.  If $p$ does not divide $\fp(\cc_j)$ for any conjugacy class $\ccj$ of $\cc$ then $p^\alpha$ divides $\fp(\ccpt)$.
\et

A well-known result, the Ito-Michler theorem, see \cite{michler}, states that the degree of every irreducible character of a finite group $G$ is coprime to a prime  $p$ if and only if a Sylow $p$-subgroup of  $G$ is both abelian and normal.

In a similar vein, for weakly-integral modular tensor categories, we obtain the following: 
\bp\label{mtc-deg-p-not div}
Let $\cc$ be a weakly-integral modular fusion category, and let $p$ be a prime divisor of $\fpcc$ such that $\fpcc=p^\al N$ with $\al\geq 1$ and $(p, N)=1$. If $p$ does not divides the Frobenius-Perron dimension  of  any of its simple objects  then $p^{\al}$ divides $|U(\cc)|$, the order of the universal grading group of $\cc$.  

Furthermore, in this case we have: $$\cc\simeq \cd\boxtimes \cd'$$  where $\cd$ is a pointed modular tensor category of dimension $p^{\al}$ and $\cd'$ is a modular tensor category of dimension $N$.
\ep

As a direct consequence of the above result we obtain the following:
\bt\label{dec-p}
Every weakly integral modular category $\cc$ decomposes as a  Deligne product $\cc\simeq \cd\boxtimes \ce$ of two weakly-integral modular categories where:
\bne
\item 
$\cd$ is pointed.
\item
$\fpcd$ and $\fp(\ce)$ are coprime integers. 
\item 
$\ce$ has the property that, for any prime divisor $p\mid \fp(\ce)$, there exists a simple object in $\ce$ whose Frobenius-Perron dimension is divisible by $p$. 
\ene
\et
For braided nilpotent fusion categories, we establish the following divisibility property:
\bt\label{br-nil}
In any braided nilpotent fusion category $\cc$, the following holds:
$$
\fp(\ccj)\mid \frac{\fpcc}{\fpccpt}.
$$
for any conjugacy class $\ccj$ of $\cc$. Moreover, the set of prime numbers dividing $\fpcc/\fpccpt$ coincides with the set of primes numbers dividing at least one element of the set  $\{\fp(\ccj)\}_{j\in \mtcj}$.
\et

The structure of this paper is as follows: In Section \ref{prelim}, we review the basic concepts of pivotal fusion categories from \cite{scalg} that are essential for the rest of the paper. Section \ref{bantay} presents several results concerning the characters of Grothendieck rings in fusion categories, which are crucial for proving Theorem \ref{br-nil}. These results may also hold independent theoretical interest. In Section \ref{dual-case}, we establish Theorem \ref{gen-dual} and its related results, Proposition \ref{mtc-deg-p-not div} and Theorem \ref{dec-p}. Finally, in Section \ref{nilpotent}, we prove Theorem \ref{br-nil} along with some other consequences.
\subsection*{Acknowledgements} 
The author thanks Sebasti\'en Palcoux  for insightful discussions and for carefully reviewing a preliminary version of this manuscript.

\section{Preliminaries}\label{prelim}
In this paper, we work over the base field $\mathbb{C}$. For a comprehensive  overview of fusion categories, we direct the reader to the monograph \cite{EGNO15}. 

Let $\cc$ be a fusion category. We denote the maximal pointed fusion subcategory of $\cc$ by $\ccpt$, and let $\irr(\cc)$ be a set of representatives for the isomorphism classes of simple objects in $\cc$. The size of $\irr(\cc)$ is referred to as the \emph{rank of $\cc$} and is denoted by $\rank(\cc)$.

A \emph{pivotal structure} in a rigid monoidal category $\cc$ is an isomorphism $j:\id_\cc\to ()^{**}$ between monoidal functors. A \emph{pivotal monoidal category} refers to a rigid monoidal category equipped with such a pivotal structure. Recall from \cite{mug1} that a pivotal fusion category is \emph{spherical} if and only if $\dim(V ) = \dim(V^*)$ for all simple objects $V$, where $\dim$ denotes the \emph{quantum dimension} of $\cc$. A fusion category $\cc$ is called \emph{pseudo-unitary} if its Frobenius–Perron dimension $\fp(\cc)$ equals its categorical dimension $\dim(\cc)$. In this case, according to \cite[Proposition 8.23]{eno-annals}, the category $\cc$ possesses a unique spherical structure, known as the \emph{canonical spherical structure}, under which the categorical dimensions of all simple objects are positive. With this spherical structure, the categorical dimension of any object $X$ in $\cc$ matches its Frobenius–Perron dimension, i.e., $\fp(X) = \dim(X)$ for every object $X$.

Recall that a fusion category is called \emph{weakly integral} if its Frobenius-Perron is an integer. Moreover, the fusion category is called \emph{integral} if the Frobenius-Perron dimension of each object is an integer.  Any weakly integral fusion category is  inherently pseudo-unitary and possesses a canonical spherical structure, under which the categorical dimensions match the Frobenius–Perron dimensions \cite[Proposition 8.24]{eno-annals}.

For any fusion category $\cc$, the monoidal center (or the Drinfeld center) of $\cc$ is a braided fusion category $\czcc$ equipped with a forgetful functor $F:\czcc\ra \cc$, see \cite{EGNO15}.

The forgetful functor $F$ also admits a right adjoint functor $R:\cc \ra \czcc$  and $Z :=FR:\cc \ra \cc$ is a Hopf comonad defined as an end.  Indeed, following \cite[Section 2.6]{scalg} one has that 
\beq
Z(V)\simeq \int_{X\in \cc}X\ot V\ot X^*
\eeq
It is known that $A:=Z(\unu)$ (the image of the unit object) has the structure of  a (central) commutative algebra in $\cz(\cc)$, and it is called \emph{the adjoint algebra} of $\cc$.

The vector space $\cecc:= \hm_{\C}(\unu, A) $ is referred to as the {\it the set of central elements}. The vector space  $\cfcc:=\hm_\cc(A,\unu)$ is referred to as the space of class functions of $\cc$. Both  $\cfcc$ and $\cecc$ have a semisimple $\comp$-algebra structure, as described in \cite{scalg}. Additionally, $\cecc$ is a commutative algebra, which means it decomposes as a direct product of copies of the field $\comp$.

Since $R:\cc \ra \czcc$ is a right adjoint to the forgetful functor $F:\czcc \ra \cc$,   \cite[Theorem 3.8]{scalg} establishes that this adjunction induces an isomorphism of $\comp$-algebras:
\beq\label{adjisom}
\cfcc \xra{\cong} \mtr{End}_{\czcc}(R(\unu)),\;\; \ch\mapsto Z (\ch)\circ \delta_\unu.
\eeq
where $\delta_\unu$ is the comultiplication of the Hopf comonad $Z$ associated at the unit $\unu$ of the fusion category $\cc$. 
\bn{exmp}
If $H$ is a semisimple Hopf algebra, then, for $\cc=\rep(H)$, as shown in \cite[Subsect. 3.7]{scalg}, it follows that $\cfcc\simeq C(H)$, the \emph{character algebra} of $H$. Recall that $C(H)$ consists of all linear trace maps $f:H\ra \comp$.
\end{exmp}
\subsection{The bijection between central primitive idempotents and conjugacy classes} \label{bij-conj} Suppose that $\cc$ is a pivotal fusion category with a commutative Grothendieck ring $K(\cc)$. Since the monoidal center $\czcc$ is also a fusion category, we can write $R(\unu)=\bigoplus_{i=0}^m\mathcal C_i$ as a direct sum of simple objects in $\czcc$. Note that any two distinct   conjugacy classes are not isomorphic to each other, since from the isomorphism of Equation \eqref{adjisom} the endomorphism ring of $R(\unu)$ is a commutative $\comp$-algebra. Thus  $\mathcal C^{0},\dots, \mathcal C^{m}$ represent the conjugacy classes of $\C$.  Since the unit object $\unu_{\czcc }$ is always a sub-object of $R(\unu)$, we can assume $\mathcal C^{0} = \unu_{\czcc }$.
In \cite[Section 5]{scalg}, Shimizu also introduced a canonical pairing between $\cfcc$ and $\cecc$ with values in $\comp$. This pairing is defined by the condition:
\beq\label{pairing-def}
f\circ a=\langle f, a\rangle=\id_\unu
\eeq
for all $f \in \cfcc$ and $a\in \cecc$.

Let $\tilde F_0, \tilde F_1, \dots, \tilde F_m\in \enx_{\czcc}(R(\unu))$ be the canonical projections onto each of these conjugacy classes. Let also $F_0, F_1, \dots F_m$ be also the corresponding primitive idempotents of $\cfcc$ under the canonical adjunction isomorphism $\cfcc \simeq \enx_{\czcc}(R(\unu))$ from Equation \eqref{adjisom}. This establishes a canonical  bijection between the (central) primitive idempotents of $\cfcc$ and the conjugacy classes of $\cc$.

Let $\mtcj:=\{0,1,\dots ,m\}$ represent the set of indices corresponding to the central primitive idempotents of 
$\cfcc$. To each $F_j$ with $j\in \mtc J$, the bijection described above assigns a conjugacy class $\ccj$.
\subsection{On the definition of conjugacy classes}
In a pivotal fusion category, Shimizu defines the elements $C_j := \mtc{F}_\lam^{-1}(F_j) \in \cecc$ as the \emph{conjugacy class sums} corresponding to the conjugacy class $\mathcal{C}^j$. Here, $\lambda \in \cf(\cc)$ is a cointegral of $\cc$ such that $\langle \lambda, u \rangle = 1$ (see \cite[Section 5]{scalg}). Here   $u:\unu\ra A$ is the unit of the central algebra $A$. Additionally, the \emph{Fourier transform} of $\cc$ associated with $\lambda$ is defined as the linear map:
\beq
\mtc F_{\lambda}:\cecc\ra \cfcc\;\;\text{given by}\;\;a \mapsto \lambda \lh \mtc S(a)
\eeq
where $\mtc{S}$ denotes the antipode.
\begin{exmp}\label{conj-class-vectg}
If $\cc = \vect^{\omega}_G$ is a pointed fusion category, then each conjugacy class sum of $G$ corresponds to a representation of the group algebra $\comp[G]$. Furthermore, the space of central elements $\cecc$ is isomorphic to $K(G)$, the Grothendieck group of $\comp[G]$. In particular, if $G$ is abelian, the conjugacy classes of $\cc$ can be identified with elements of $G$, since in this case, the dual group $\widehat{G}$ is isomorphic to $G$.
\end{exmp}
\subsection{The integral of a fusion category}
Let $\cc$ be a fusion category, and let $A=Z(\unu)$ be its adjoint algebra as defined earlier. Following \cite[Definition 5.7]{scalg}, an \emph{ integral} in $\C$ is a morphism $\Lambda: 1 \ra A$ in $\C$ that satisfies
\beqn
 m \circ (\id_{A}\otimes \Lambda)=\eps_{\unu}\ot \Lam.
\eeqn
where $m:A\otimes A\ra \unu$ is the multiplication of the adjoint algebra $A$. Here, $\eps_\unu$ is counit $\eps_V:Z(V)\ra V $ specialized at $V=\unu$, the unit of the fusion category $\cc$.

It is well known that the integral of a finite fusion category is unique up to a scalar.  Let $\Lam\in \cecc$ be the idempotent integral of $\cc$ and define $\widehat{\Lam}:=\dim(\cc)\Lam$. 
 
Let $\{V_0, \dots, V_m\}$ be a complete set of representatives of the isomorphism classes of simple objects in $\cc$. Shimizu, in \cite[Section 6]{scalg}, provided a construction of the (central) primitive elements
$E_i\in \cecc$,  corresponding to the irreducible characters $\ch_i:=\ch_{V_i}$ associated to each simple object $V_i$ of $\cc$. This correspondence is defined by the condition: $$\lag\ch_i, E_j\rag=\delta_{i, j},$$ for all $i,j$.
 
Note that \cite[Lemma 4.1]{ccc-march}, together with Equation (4.7) from the same paper, provides an alternative pair of dual bases for the canonical pairing. More precisely, it shows  that $\{F_j\}$ and $\frac{C_j}{\dim(\ccj)}$ form another pair of dual bases expressed as:
\beq\label{dbase}
< F_j, \frac{C_{j'}}{\dim(\cc^{\jp})}>=\delta_{j,j'}
\eeq
for all $j, j'\in \mtc J$.

Without loss of generality, we can assume that $V_0=\unu$, the unit of the fusion category $\cc$. This gives $\ch_0=\epsu$, and the associated primitive idempotent is $E_0=\Lam_\cc$, see \cite[Lemma 6.1]{scalg}.

The following lemma is straightforward, but for reader's convenience, a proof is provided.
\bl
Let $\cc$ be a pivotal fusion category with a commutative Grothendieck ring. We have the following formula for the integral element of $\cc$:
\beq\label{integrform}
\widehat{\Lam}_{\cc}=\sum_{j\in \mtcj} C_j.
\eeq 
\el
\bpf
We will show that both sides of the  Equation \eqref{integrform} have the same image under $\mtf$. For the right hand side, by the definition of conjugacy class sums $C_j$,  we have  $\mtc F_\lam(\sum_{j\in \mtcj} C_j)=\sum_{j\in \mtcj} F_j=\eps_\unu$. For the left hand side, according to \cite[Equation (6.10)]{scalg}, we have 
\beq\label{610}
\mtf(E_{i^*})=\frac{d_i}{\dimcc}\ch_i
\eeq
for any $0\leq i\leq m$. In the case $i=0$, one has $\ch_0=\epsu$ and $E_0=\blam_\cc$. Since $d_0=1$ it follows that $\mtf(E_{0})=\epsu.$
\epf

\subsection{Definition of the support $\mtc J_{\cd}$}
Let $\cc$ be a fusion category whose Grothendieck ring is commutative. For any fusion subcategory $\cd \subseteq \cc$, recall the notion of its support from \cite[Subsection 4.2]{ccc-march}. It is well known that there is an inclusion of $\comp$-algebras $\iota: \cf(\cd) \hookrightarrow \cf(\cc)$. Consequently, there exists a subset $\mtc{J}^\cc_\cd \subseteq \mtc J^\cc:=\{0, \dots, m\}$ such that
\beq\label{idemptint}
\lam_{\cd}=\sum_{j \in \mtc J_\cd}F_j
\eeq
since $\lambda_\cd$ is an idempotent element in $\cf(\cc)$. The set $\mtc{J}_\cd$ is referred to as the \emph{support} of the fusion subcategory $\cd$.

For a fusion category $\cc$, the fusion subcategory $\ccad$ represents the trivial component in the universal grading of $\cc$ (see \cite{NG}). Alternatively, $\ccad$ is the smallest fusion subcategory of $\cc$ which contains all the objects $X\otimes X^*$, with $X\in \irr(\cc)$.

A grading of a fusion category $\cc$ by a group $G$
is defined by a map $\deg : \irr(\cc) \ra G$ with the property that  for any  simple objects $X, Y, Z $ of $\cc$,  if $X\otimes Y$ contains $Z$ then $\deg(Z) =\deg(X) \deg(Y)$.  The term \emph{"grading"}  also refers to the  decomposition $\cc=\bigoplus_{g\in G} \cc_g$, where $\cc_g $ is the full additive subcategory generated by the simple objects of degree $g \in G$.

The subcategory $\cc_1$ corresponding to $g =1$,  is called the \emph{trivial component of the grading} and is a fusion subcategory. A grading is \emph{trivial} if $\cc_1=\cc$, and  \emph{faithful} if the map $\deg : \irr(\cc) \ra G$ is surjective. The set of all possible gradings
of $\cc$ by $G$ depends functorially on $G$, leading to the concept of a \emph{universal grading}. As shown in \cite[Section 3.2]{NG}, such a grading always exists, and the associated group is denoted by $U(\cc)$. The following  results are proven in \cite[Section 3.2]{NG}:
(see \cite[(Lemma 3.4, Theorem 3.5, and Proposition 3.9)]{NG}.
\bne
\item 
Any  faithful grading of $\cc$ by a group $G$ is determined by a surjective group homomorphism $\pi :U(\cc)\ra G$.
\item 
The universal grading $\deg : \irr(\cc) \ra U(\cc)$ is faithful.
\item 
The trivial component of the universal grading equals $\ccad$. 
\item
Any fusion subcategory $\cd\subseteq \cc$ that contains $\ccad$ can be written as form $\cd=\bigoplus_{g\in G}\cc_g$ for some subgroup $G\leq U(\cc)$.
\ene
It follows that $\fp(\cc)=|U(\cc)|\fpccad$.

Let $\cc$ be a fusion category. Following \cite{NG}, we define a sequence of subcategories: $\cc^{(0)}=\cc$, $\cc^{(1)}=\ccad$, and  for all $n\geq 1$,  $\cc^{(n)}=\cc^{(n-1)}_{ad}$.

The descending sequence 
$$
\cc=\cc^{(0)}\supseteq  \cc^{(1)}\supseteq \dots \supseteq \cc^{(n)} \supseteq \dots
$$
is called the \emph{upper central series}. A fusion category $\cc$ is said to be \emph{nilpotent} if its upper central eventually reaches  $\vect$ (the category of finite dimensional $\comp$-vector spaces), meaning that  $\cc^{(n)}= \vect$ for some $n$.  The smallest such
$n$ is called the \emph{nilpotency class} of $\cc$.
\section{Restriction of characters of Grothendieck rings}\label{bantay}
Throughout this section, let $\cc$ be a fusion category with a commutative Grothendieck ring, $K(\cc)$. Then $\cfcc$ is a commutative semisimple algebra, and $\cfcc\simeq K(\cc)$, as shown in \cite[Example 4.4]{scalg}.  
As above, Let $\mtc J^{\cc}:= \{0, \dots, m\}$ represent the set of indices corresponding to the conjugacy classes of $\cc$. It is clear that the rank of $\cc$ is $\rank(\cc) = |\mtc J^\cc| = m + 1$.

Suppose $\cd$ is a fusion subcategory of $\cc$. Note that the Grothendieck ring $K(\cd)$ is also commutative, and, as in the case above, $K(\cd)\simeq \cf(\cd)$ as algebras.  We let $\mtc J^{\cd}:= \{0, \dots, m'\}$ represent the set of indices corresponding to the conjugacy classes of $\cd$. As above, it is clear that the rank of $\cd$ is $\rank(\cd) = |\mtc J^\cd| = m + 1$.

For any object $V \in \cc$, following \cite[Sect. 4.3]{scalg}, consider the end
$$
\barz(V):\cd^{\op}\times \cd\ra \cc,\;\; \barz(V):=\int_{X\in \cd}  X\ot V\ot X^*.
$$ 
Let $\bar{\pi}_{V;\; X}:\barz(V) \ra X\ot V\ot X^*$ denote the universal dinatural maps that define this end. By the universal property of $\barz$, for any $V\in \cc$ there  exists a unique canonical map in $\cc$:
\beq
Z(V)\xra{q_V}\barz(V)
\eeq
such that $\bar{\pi}_{V;\; X}\circ q_V={\pi}_{V;\; X}$ for every  object $X$ in $\cd$.

In \cite[Appendix]{ccc-march} we demonstrated  that $q$ is a morphism of Hopf comonads. The map $q_\unu:Z(\unu)\ra \barzu$ induces two associated maps:
$$
\pif=q_\unu^*:\cfcd \ra \cfcc, \;\ch\sent \ch \circ q_\unu,
$$
$$
\pie={q_{\unu}}_*:\cecc \ra \cecd, \;z\sent q_\unu \circ z.
$$ 

Furthermore, by \cite[Lemma 7.7]{ccc-march} it is known that $q_\unu^*$ is a monomorphism and ${q_\unu}_*$ is an epimorphism. Both maps are, in fact, $\comp$-algebra homomorphisms.

\bn{exmp}
If $\pi:H\ra K$ is a surjective Hopf algebra homomorphism, then clearly we have $\rep(K)\subseteq \rep(H)$. In this case, $\pi_e=\pi\big|_{Z(H)}:Z(H)\ra Z(K)$ is the restriction of $\pi$ to the center $Z(H)$ of $H$. Furthermore, $\pi_f$ corresponds to the  canonical inclusion  of the character algebra $C(K)$ of $K$ into that $H$, i.e, $C(K)\hookrightarrow C(H)$.
\end{exmp}

Let, $\{F_j\}_{0 \leq j \leq m}$ be  the (central) primitive idempotents of $\cfcc$, and let $\{\cc^j\}_{0\leq j \leq m}$ be their associated conjugacy classes. For each (central) primitive  idempotent $F_j\in \cfcc$, let $\muj:\cfcc\ra \comp$ denote its associated linear character. Similarly, let $\{G_t\}_{0 \leq t \leq m'}$ be the (central) primitive idempotents of  $\cf(\cd)$, with  $\{\cd^t\}_{0\leq t \leq m'}$ being  their associated conjugacy classes.  Additionally, let $\nu_t: \cfcd \to \comp$ denote the characters associated with these idempotents.

Using the map $\pif$ we can identify $\cfcd$ as subring of  $\cfcc$.

 From the proof of \cite[Theorem 6.12]{b-galois}, for all $\al \in \cfcc$, we have 
\beq\label{dual-wkc}
\muj(\al)=\langle \al, \frac{C_j}{\dim(\ccj)}\rangle
\eeq

We say that two characters $\muj, \mu_{j'}:\cfcc\ra \comp$ are equivalent or \emph{in the same class},  (denoted $\mu_j \simeq \mu_{j'}$) with respect to $\cd$ if their restrictions to $\cf(\cd)$ are the same.

The set of equivalence classes is denoted by $\clcrd$ and it is referred to as \emph{the set of  classes} of $\cc$ with respect to $\cd$. This set is   in bijection with the set  of the characters of $\nu_t:\cfcd\ra\comp$ with $0\leq t\leq m'$.  
 
Moreover, an element $C\in  \clcrd$ consists of all the indices $j$ such that the restriction $\mu_j\dw_{\cfcd}$ corresponds to the same character in $\cfcd$, denoted as $\ro_C:=\nu_s$ for some $s\leq m'$. Thus, we use $\ro_C$ to represent the restriction to  $\cfcd$ of all characters $\mu_j$ belonging to some class $C\in \clcrd$. 

As noted above, the set of character classes of $\cc$ relative to $\cd$ corresponds bijectively to the characters of $\cfcd$. Consequently,  the central primitive idempotents of $\cfcd$ are also indexed by the classes $C\in \clcrd$. We denote the central primitive idempotent associated with a class $C \in \clcrd$ by $G_C \in \cf(\cd)$. Thus, the elements $G_C$ for $C \in \clcrd$ form a linear basis for $\cf(\cd)$, consisting of the complete set of (central) primitive idempotents of the commutative $\comp$-algebra $\cf(\cd)$. For the remainder of this section, we identify the set of indices $\mtc J^\cd$ of the central primitive idempotents of $\cf(\cd)$ with the elements of $\clcrd$. For any primitive idempotent $G_C$ of $\cfcd$, as described in Subsection \ref{bij-conj}, let  $\mtr{D}_C\in \cecd$ denote  the corresponding conjugacy class sum of $\cd$ associated with $G_C$. 
\bp \label{nec-p} 
Let $\cc$ be a pivotal fusion category with a commutative Grothendieck ring, $K(\cc)$, and let $\cd$ be a fusion subcategory. With the above notation, we have
\beq\label{pizgm}
\pi_f(G_C)=\sum_{j \in C}F_j
\eeq
for any class $C$.
\ep
\bpf
Since $\pi_f$ is an algebra homomorphism and $\pi_f(G_s)$ is idempotent, it can be expressed as a sum of primitive idempotents in $\cfcc$. Suppose that $\piz(G_C)=\sum_{j \in \mtc B_C}F_j$ for some subset $\mtc B_C\subseteq \{0,1,\dots, m\}$. This implies that $\mu_j(G_m)=\mu_{j'}(G_m)=1$ for any $j,j' \in \mtc B_C$, indicating that $j$ and $j'$ belong to the same character class, denoted as $C$. Indeed, we have  $\mu_j(G_t)=\mu_{j'}(G_t)=0$ for any $t\neq s$, which leads to $\muj|_{\cfcd}=\mu_{j'}|_{\cfcd}$. 

Conversely, if $j'\in C$ then $\mujp\dw_\cfcd=\muj\dw_\cfcd$. This implies that $\mujp(G_C)=1$, and therefore $\jp\in \mtc B_C$.
\epf

Recall the canonical pairing between $\cfcc$ and $\cecc$, which we denote as  $\langle,\;\rangle_\cc$. For the fusion subcategory $\cd\subseteq \cc$, there is an analogous pairing, written here as  $\langle,\;\rangle_\cd$. The proof of the following lemma follows directly from the definitions of the two canonical pairings.
\bl \label{l1}
The two canonical pairings on the fusion categories $\cc$ and $\cd$ are compatible, satisfying the relation:
\beq\label{comp1}
<\pif(\ch), \;z>_\cc=<\al, \piz(z)>_\cd
\eeq
for all $\ch \in \cfcd$ and $z\in \cecc$.
\el
\subsection{On the projection of conjugacy class sums}
\bt\label{projclass} Let $\cc$ be a pivotal fusion category with a commutative  Grothendieck ring,  and $\cd$ be a fusion subcategory of $\cc$. Let $j\in \mtcj$ be any index  and suppose that $j\in C$ for some $C\in \clcrd$. Then:
\beq\label{piz-form}
\piz(\lbarcj)=\frac{\mtr{D}_C}{\dim(\cd^C)}
\eeq
\et
\bpf  
We can write
$$
\piz(\lbarcj)=\sum_{C_1 \in \clcrd}\al_{jC_1}\mtr{D}_{C_1}
$$
for some scalars $\al_{jC_1}\in \comp$.

Let now $C_2\in \clcrd$ be an arbitrary class, and let $\ch=G_{C_2}$ and $z= \mtr C_j$. By Lemma \ref{l1} we have:
$$
<\sum_{s \in C_2} F_s, \lbarcj>_{ _\cc}=<G_{C_2}, \pi_e(\lbarcj)>_{_\cd},
$$
which in turn can be written as 
$$
<\sum_{s \in C_2} F_s, \lbarcj>_{ _\cc}=<G_{C_2}, \sum_{C_1\in \clcrd}\al_{jC_1}\mtr{D}_{C_1}>_{_\cd}.
$$

The right hand side of the last above equation, using the dual bases form Equation \ref{dbase}  for $\cd$, can be written as: 
$$
<G_{C_2}, \sum_{C_1 \in \clcrd}\al_{jC_1}\mtr{D}_{C_1}>_{ _\cd}=\al_{jC_2}\dim(\cd^{C_2}).$$

On the other hand, using  also the dual bases form Equation \ref{dbase}, for $\cc$ this time, the left hand side equals $\delta_{j\in C_2}$.

Thus, we find that $\al_{jC_2}\dim(\cd^{C_2})=\delta_{j\in C_2}$ for all character  classes $C_2\in \clcrd$.  Since $j\in C$, this gives that $\alpha_{jC} = \frac{1}{\dim(\cd^C)}$ and $\al_{jC_2}=0$ for all $C_2\neq C$. This establishes Equation \eqref{piz-form}.
\epf

For a fusion category $\cc$ we define the natural coend projection $\eps_\cc:A\ra \unu\in \cfcc$ as the map  $A\xra{\pi_{\unu, \unu}}\unu \ot \unu^*\simeq \unu$. This projection is referred to as  $\epsu$ in \cite{scalg}, but we use the notation $\eps_\cc$ to emphasize the category $\cc$. 

Recall that in the previous section, we defined $\ch_0\in \cfcc$ as the irreducible character corresponding to  the tensor unit $\unu_\cc$ of $\cc$. Without loss of generality, we can also assume that  $\al_0\in \cfcd$ is the character associated to the unit object $\unu_\cd$ of $\cd$. It is clear that from its definition that $\pif(\al_0)=\ch_0$.

\bl\label{pioflam}
Let $\cd\subseteq \cc$ be a fusion subcategory of a pivotal fusion category $\cc$. Under the above notation, we have:
\beq\label{epsqu}
\eps_\cd\circ \qu=\eps_\cc
\eeq
and
\beq\label{pizlam}
\piz(\blam_\cc)=\blam_\cd
\eeq
\el
\bpf
From the naturality of $\qu$, we can conclude that $\eps_\cd\circ \qu=\eps_\cc$, since $\eps_\cc:A\ra \unu$ serves as the universal end projection $A\xra{\pi_{\unu, \unu}}\unu \ot \unu^*\simeq \unu$.

For the second assertion, we first show that $\piz(\blam_\cc)\neq 0$. According to \cite[Lemma 6.1]{scalg},  the corresponding idempotent of $\ch_0$ is $E_0=\blam_\cc$. By applying Lemma \ref{l1}, we have:
$$
1=<\pif(\al_0), E_0>=<\al_0, \piz(\blam_\cc)>,
$$ 
which implies that $\piz(\blam_\cc)\neq 0$. Moreover, for any $z \in \cecc$, we find:
$$
\eps_\cc(z)\piz(\blam_\cc)=\piz(z\blam_\cc)=\piz(z)\piz(\blam_\cc)
$$

Note that $\eps_\cd(\pi(z))=\eps_\cc(z)$, and since $\piz$ is surjective, it follows that:
$$
w\piz(\blam_\cc)=\eps_\cd(w)\piz(\blam_\cc)
$$
for any $w\in \cecd$. Given that $\piz(\blam_\cc)\neq 0$ it follows  $\piz(\blam_\cc)=\blam_\cd$.
\epf
\section{On the dimensions of conjugacy classes}\label{dual-case}
For any weakly integral fusion category, whose Grothendieck group $K(\cc)$ is commutative, the following Harada type identity holds, as established in \cite[Theorem 1]{harada}: 
\beq\label{harada}
(\prod_{j=0}^m \frac{C_j}{\fp(\cc^j)})^2=\frac{\fp(\ccpt)}{\fpcc}\big(\sum_{j\in \mtc J_{\ccpt}}C_j\big)
\eeq
Recall that the support of $\cd$ in $\cc$ is the subset of indices of  $\mtc J_{\cd}\subseteq \mtc J$ such that $\lam_{\cd}=\sum_{j \in \mtc J_\cd}F_j$.

Since the elements $\{C_j\}$ form a linear basis for $\cecc$, we can express the product of two basis elements as:
\beq\label{stc-def}
C_iC_j=\sumktom c^k_{ij}C_k.
\eeq
where $c^k_{ij}\in \mathbb C$ are \emph{the structure constants} of $\cc$. It is known that these constants are in fact integers in the case $\cc=\rep(G)$ for  a finite group $G$. Burnside provided a formula relating these structure constants to the values of irreducible characters, see \cite[pp. 45]{is}.  A similar result for the structure constants of a pivotal fusion category with a commutative Grothendieck ring was presented in \cite[Theorem 1.1]{b-blms}. 

A braided fusion category with a spherical structure is known as a \emph{premodular} fusion category. The structure constants for premodular categories were investigated in \cite{b-blms}.

In the case of premodular categories, a specific formula for the structure constants $c^k_{ij}$ was provided in the proof of \cite[Theorem 1.2]{b-blms}. It was demonstrated that: 
\beq
\mtr C_i\mtr C_j=\sum_{k=1}^m\bigg(\sum_{v\in \ca_k}\frac{{\wdht N}^v_{\sg(i)\sg(j)}\widehat{d}_v}{\dim(\cc^k)}\bigg)\mtr C_k
\eeq
which  leads to the following expression for the structure constants:
\beq\label{str-form}
c^k_{ij}=\sum_{v\in \ca_k}\frac{{\wdht N}^v_{\sg(i)\sg(j)}\widehat{d}_v}{\dim(\cc^k)}.
\eeq

By \cite[Prop 8.23]{eno-annals}, since every weakly integral braided fusion category is premodular, the results of \cite[Theorem 1.2]{b-blms} apply in this context. Therefore from Equation \eqref{str-form} we can conclude that:
\beq\label{more-rat}
\dim(\cc^k)c^k_{ij}\in \mathbb Z_{\geq 0}
\eeq
for any structure constant $c^k_{ij}$ in a weakly integral braided  fusion category. Note the that the fusion rules ${\wdht N}^v_{\sg(i)\sg(j)}$ are integers, and since $\cc$ is weakly integral, the quantum dimensions $\widehat{d}_v$ are also integers. For further details, see the proof of \cite[Theorem 1.2]{b-blms}. 
\subsection{Proof of Theorem \ref{gen-dual}}
\bpf 
Let $\wtp:=\prod_{j=0}^m\fp(\cc^j)$. Since $\cc$ is weakly integral,  by \cite[Proposition 8.27]{eno-annals} implies that the Frobenius–Perron
dimensions of simple objects in $\ccad$ are integers. Given that the conjugacy classes $\ccj$ are direct sums of simple objects of the adjoint subcategory, it follows that $\fp(\ccj)$ are also integers and therefore $\wtp\in \mathbb Z$.

From the proof of \cite[Theorem 1.2]{b-blms}, we know that $\dim(\cc^k)c^k_{ij}\in \mathbb Z_{\geq 0}$ for any premodular fusion category. Therefore, multiplying equation \eqref{stc-def} by $\wtp$, we have:
$$
\wtp C_iC_j=\sumktom (\wtp c^k_{ij})C_k
$$
This shows that the coefficients $\wtp c^k_{ij}\in \mathbb Z$. By induction, for any $r$-class sums,  the product
$\wtp^{\;r-1}C_{i_1}C_{i_2}\dots C_{i_r}$  yields integral coefficients for each $C_k$ when expressed as a linear combination of conjugacy class sums. 

Next, by multiplying by $\wtp^{\;2m+1}$ in Equation \eqref{harada}, we obtain:
\beq\label{harada-braided}
\wtp^{\;2m+1}(\prod_{j=0}^m {C_j})^2=\frac{\wtp^{\;2m+3} \fp(\ccpt)}{\fpcc}\big(\sum_{j\in \mtc J_{\ccpt}}C_j\big)
\eeq
Thus, in the product $\wtp^{\;2m+1}(\prod_{j=0}^m {C_j})^2$, all $C_j$ with $j\in \mtc J_\ccpt$ have a integral coefficient. Thus, $\frac{\wtp^{\;2m+3}\fp(\ccpt)}{\fpcc}\in \mathbb Z_{>0}$ and the proof is complete since $\rank(\cc)=m+1$.
\epf
\subsection{Fusion categories with rational structure constants}
\bt\label{D}
Let $\cc$ be a weakly integral fusion category, and suppose that all its structure constants $c^k_{ij}$ are rational numbers. Let $D\in \mathbb Z$ be such that $Dc^k_{ij}\in \mathbb Z$ for all $i, j, k$. Then:
\bne
\item 
$\large{\frac{D^{2m+1} \wtp^{\;2}\fp(\ccpt)}{\fpcc}\in \mathbb Z}$.
\item $d_i\mid D\;\fp(\cc)$
\item
$d_i\mid D^{2m+2} \wtp^{\;2}\fp(\ccpt)$.
\ene
\et
\bpf
The proof of the first item follows the same method as the previous theorem, but this time multiplying  Harada's identity by $D^{2m+1}$. This ensures an integral coefficient for $C_0$ on the left hand side of the equation. The second item is a direct consequence of  \cite[Theorem 4.7]{b-palcoux}. The third item follows directly from the first two divisibility results.
\epf
\bl
Let $\cc$ be a weakly integral fusion category, and suppose that all its structure constants $c^k_{ij}$ are rational numbers. Let $\cd$ be a fusion subcategory of $\cd$. Then $\cd$ has also rational structure constants.
\el
\bpf
Recall the map $\pi_e:\cecc\ra\cecd$ from Section \ref{bantay}. Let $\cd_E$ and $\cd_F$ be two conjugacy classes of $\cd$ with $E, F\in \clcrd$. Since $\pi_e$ is a surjective map,  as in Equation \eqref{piz-form}, there are  indices  $s\in E,\; t\in F$ such that 
$
\pi_e(\mtr{C}_{t})=\frac{\dim(\cc^{t})}{\dim(\cd^t)}\mtr{D}_E
$
and 
$
\pi_e(\mtr{C}_{s})=\frac{\dim(\cc^{s})}{\dim(\cd^s)}\mtr{D}_F.
$  
Then
$$
\mtr{D}_E\mtr{D}_F=\frac{\dim(\cd^E)}{\dim(\cc^{t})}\frac{\dim(\cd^F)}{\dim(\cc^{s})}\pi_e(\mtr{C}_{t})\pi_e(\mtr{C}_{s})
$$
On the other hand, since $\pi_e$ is an algebra homomorphism one has 
\begin{eqnarray*}
\pi_e(\mtr{C}_{t})\pi_e(\mtr{C}_{s}) &=& \pi_e(\mtr{C}_{t}\mtr{C}_{s})=\pi_e(\sumktom  c^k_{ts}\mtr{C}_k)\\ &= &\sumktom c^k_{ts}\frac{\dim(\cc^{k})}{\dim(\cd^{\pi(k)})} \mtr{D}_{\pi(k)}\\ &= &\sum_{C\in \clcrd}\big(\sum_{k\in C}c^k_{ts}\frac{\dim(\cc^{k})}{\dim(\cd^{C})}\big) \mtr{D}_{C}\\ &= & \sum_{C\in \clcrd}\big(\sum_{k\in C}c^k_{ts}\dim(\cc^{k})\big)\frac{ \mtr{D}_{C}}{\dim(\cd^{C})}
\end{eqnarray*}
Therefore 
\begin{eqnarray*}
\mtr{D}_t\mtr{D}_s &=& \frac{\dim(\cd^t)}{\dim(\cc^{t_1})} \frac{\dim(\cd^s)}{\dim(\cc^{s_1})}\sum_{C\in \clcrd}\big(\sum_{k\in C}c^k_{t_1s_1}\dim(\cc^{k})\big)\frac{ \mtr{D}_{C}}{\dim(\cd^{C})}
\end{eqnarray*}
which shows that the coefficient of each class sum $ \mtr{D}_C$ is a rational number.
\epf
\bp
Let $\cc$ be a weakly integral fusion category with a commutative Grothendieck ring and suppose that all its structure constants are integral, i.e., $c^k_{ij}\in \mathbb Z$. With above notation, one has:
\beq\label{div:fpcc-int}
\frac{\fp(\ccpt)}{\fp(\cc)}\big[\prod_{j=0}^m\fp(\ccj)\big]^{2}\in \mathbb Z.
\eeq
Moreover, if $U(\cc)$ has odd order, then we have: 
\beq\label{div:fpcc-odd-int}
\frac{\fp(\ccpt)}{\fp(\cc)}\big[\prod_{j=0}^m\fp(\ccj)\big]\in \mathbb Z.
\eeq
\ep
\bpf
If $\cc$ has integral structure constants then, we can take $D=1$ since the coefficient of $C_0$ in the left hand side of Equation \eqref{harada} is already an integer.  Therefore, in this case, there is no need  to multiply the equation by $\wtp^{\;2m+1}$, and we  directly obtain:
$$
\frac{\fp(\ccpt)}{\fp(\cc)}\big[\prod_{j=0}^m\fp(\ccj)\big]^{2}\in \mathbb Z.
$$

Using the dual version of \cite[Corollary 4.11]{b-pa}, in the case in which $U(\cc)$ has odd order, we also have the second divisibility result. 
\epf
A well-known result regarding the prime factors of the sizes of conjugacy classes in finite groups can be found, for example, in \cite[Proposition 4]{herzog}. Let $G$ be a finite group and $p$ a prime divisor of $G$. Then, $p$ does not divide  the size of any conjugacy class of $G$ if and only if the Sylow $p$-subgroup of  $G$ is contained in its center. In this case, Schur–Zassenhaus theorem, see \cite{robinson}, ensures that $G\simeq H\times K$ where $H$ is its Sylow p-subgroup and $K$ is another subgroup of $G$. In particular this shows that $p^{\al}$ divides the order of the abelian group $G/G'$, i.e. dimension of $\rep(G)_{pt}$

As a consequence of  Theorem \ref{gen-dual}, a weaker version of this result can be proven for weakly integral braided  fusion categories in Theorem \ref{ccj-not div-p}.

\subsection{Proof of Theorem \ref{ccj-not div-p}}
\bpf
Since $\frac{\wtp^{\;2m+3}\fp(\ccpt)}{\fpcc}\in \mathbb Z_{>0}$, and  $p$ does not divide $\wtp$,  we conclude that $p^{\al}\mid \fp(\ccpt)$.
\epf
For finite groups we obtain the following:
\bc
Let $G$ be a finite group of order $|G|=p^\al N$ with $\al\geq 1$ and $(p, N)=1$. If $p$ does not divide the size of any conjugacy class of $G$, then $p^{\al}$ divides $\frac{|G|}{|G'|}$. 
\ec
\bpf
Corollary \ref{ccj-not div-p} implies that $p^{\al}\mid \fp(\ccpt)=\frac{|G|}{|G'|}$.
\epf

\subsection{Proof of Proposition \ref{mtc-deg-p-not div}}
\bpf
Since $\cc$ is non-degenerate, we have $\czcc\simeq \cc\bxt \cc^{\rev}$, as shown in  \cite[Proposition 3.7]{dgno}. Thus, any simple object in $\czcc$ is of the form $V_i\bxt V_{m}$, and from the hypothesis, the prime  $p$ does not divide the Frobenius-Perron dimension of any such simple object. Given that the conjugacy classes of $\cc$ are simple objects in $\czcc$, it follows that the conditions of Theorem \ref{ccj-not div-p} are met. Consequently, $p^{\al}$ divides $\fp(\ccpt)$. Furthermore, in any modular tensor category, it is known from \cite[Theorem 6.2]{NG} that $\gcc\simeq U(\cc)$, which implies that $p^{\al}$ divides the order of $U(\cc)$.

Since $p^{\al}\mid \fp(\ccpt)$ there exists a fusion subcategory a $\cd\subseteq \ccpt$ with $\fp(\cd)=p^{\al}$. Let $\cd'$ be the M\"uger centralizer of $\cd$ in $\cc$. By \cite[Theorem 3.10]{dgno}, we know that $\fp(\cd')=N$, and thus $\cd\cap \cd'=\vect$. Furthermore, \cite[Corollary 3.5]{dgno} implies that $\cd$ is braided non-degenerate and therefore $\cc\simeq \cd\bxt \cd'$. According to \cite[Theorem 3.10]{dgno} $\cd'$ is also non-degenerate.
\epf
\br
The fact that $p\mid \fp(\ccpt)$  also follows from \cite[Equation (1.14)]{b-pa}. It is worth noting  that for finite groups,  as shown in \cite{michler}, the proof of the Ito-Michler Theorem is significantly more complex and relies on the Classification of Finite Simple Groups.
\er
\subsection{Proof of Theorem \ref{dec-p}}

\bpf
We apply Proposition \ref{mtc-deg-p-not div} repeatedly for each prime that does not divide the dimension of any simple object in $\cc$. Suppose $\fpcc = p_1^{\alpha_1} \dots p_n^{\alpha_n}$. Without loss of generality, assume that the primes $p_i$ for $1 \leq i \leq r$ do not divide the dimension of any simple object in $\cc$.

Applying Proposition \ref{mtc-deg-p-not div} with $p = p_1$, we obtain a decomposition $\cc = \cd_1 \boxtimes \ce_1$, where $\cd_1$ is a pointed modular fusion category with $\fp(\cd_1) = p_1^{\alpha_1}$. Since $\fp(\ce_1) = p_2^{\alpha_2} \dots p_n^{\alpha_n}$, we can apply Proposition \ref{mtc-deg-p-not div} again with $p = p_2$, yielding $\ce_1 = \cd_2 \boxtimes \ce_2$, where $\cd_2$ is a pointed modular fusion category with $\fp(\cd_2) = p_2^{\alpha_2}$. 
 
Continuing in this way for each prime $p_i$ with $1 \leq i \leq r$, we obtain
$$
\cc\simeq \cd_1\boxtimes \dots cd_r\boxtimes \ce_r.
$$
Thus we set $\cd:=\cd_1\boxtimes \dots \cd_r$ and $\ce:=\ce_r$.
\epf
\begin{corollary} 
For any weakly integral perfect modular fusion category $\cc$, if $p$ is a prime divisor of $\fpcc$, then there exists a simple object in $\cc$ whose squared dimension is divisible by $p$. 
\end{corollary}
\bn{exmp}\label{counter-example}
Let $\ca=\rep(D(S_3))$ be the fusion category of representations of the Drinfeld double of the symmetric group $S_3$. As described, for instance, in  \cite{CW4}, this category contains $2$ invertible objects, $4$ simple objects of dimension $2$, and $2$ simple objects of dimension $3$.  Additionally,  the $2$ invertible objects together with the $4$ simple objects of dimension $2$ form a fusion subcategory $\cc$ with Frobenius-Perron dimension $\fp(\cc)=2.1^2+4.2^2=18$.  In fact, it can be verified that $\cc=\ca_{ad}$. 

For $p=3$ we have $p\mid \fp(\cc)$ (with $\al=2$). Moreover, $p$ does not divide the dimension of any simple object in $\cc$,  and $p\nmid \fp(\cc_{\mtr{pt}})=2$. Therefore, in Proposition \ref{mtc-deg-p-not div},  the assumption that $\cc$ is non-degenerate cannot be omitted. From Theorem \ref{ccj-not div-p}, it also follows that $\cc$ has at least one conjugacy class with a dimension divisible by $3$.
\end{exmp}
The example above, $\cc = \rep(D(S_3))_{\ad}$, also illustrates that, in general, the set of prime divisors of the Frobenius-Perron dimensions of simple objects does not always coincide with the set of prime divisors of the Frobenius-Perron dimensions of the conjugacy classes.
\section{Proof of Theorem \ref{br-nil}}\label{nilpotent}
We would like to thank A. Moreto for suggesting the proof of the following result:
\bp\label{R1-nilpotent}
For any weakly integral fusion category with a commutative Grothendieck ring  $K(\cc)$, the following inequality holds:
\beq\label{ineq}
\fp(\ccj)\leq \frac{\fpcc}{\fpccpt},
\eeq
for any conjugacy class $\ccj$ of $\cc$.
\ep

\bpf
Since $\ccpt\subseteq \cc$, by \cite[Lemma 7.5, item 1)]{ccc-march} there exists a surjective homomorphism of $\comp$-algebras:
$$
\pi_e: \cecc \ra \mtr{CE}(\ccpt),
$$
as discussed in Section \ref{bantay}.

Let $\ccpt:=\vect_G^\omega$ for some abelian group $G$ and some $\omega\in H^3(G, \comp^*)$. It follows from Example \ref{conj-class-vectg}, we have $\mtr{CE}(\ccpt)= \comp[G]$. Note that the conjugacy classes of $\vect_G^\omega$ correspond to the elements of $G$. 

According to Theorem \ref{projclass} for any conjugacy class 
$\ccj$ of $\cc$, there exists a group element $g_j\in G$ such that: 
$$
\pi_e(C_j)= \dim(\ccj) g_j.
$$
Since $\cc$ is weakly integral, it is a pseudo-unitary fusion category  implying that  $\dim(\cc^j)=\fp(\ccj)>0$ for any $j\in \mtc J$.

Additionally,  by Lemma \ref{pioflam}, we have:
\beq\label{wtlam}
\pi_e( \widehat{\Lam}_\cc) =\frac{\fpcc}{\fp(\ccpt)} \widehat{\Lam}_\ccpt.
\eeq
Since $\mtr{CE}(\ccpt)= \comp[G]$, it follows that $\widehat{\Lam}_\ccpt=\sum_{g \in G} g$. Recall that, according to \cite[Proposition 4.1]{ccc-march} (applied for $\cd=\cc$) we have: 
$$
\widehat{\Lam}_\cc=\sum_{j\in \mtc J} C_j.
$$
By the same \cite[Proposition 4.1]{ccc-march} applied for $\cd=\ccpt$ , we obtain:
$$ 
\widehat{\Lam}_\ccpt=\sum_{j \in \mtc J_{\ccpt} }C_j.
$$
Let $\mtc B_j:=\{l \in \mtc J\;|\; \pi_e(C_l)=\fp(\cc^l)g_j\}$.
Counting the multiplicity of $g_j$ in the left hand side of Equation \eqref{wtlam} we find:
$$
\sum_{l \in {\mtc B}_j} \fp(\cc^l) = \frac{\fpcc}{\fpccpt}.
$$
Since $\fp(\ccl)>0$ for any $l\in \mtc J$, we conclude that:
$$\fp(\ccl)\leq  \frac{\fpcc}{\fpccpt}.$$
\epf
\br
The above proof shows that the possible equality $\dim(\cc^j) = \frac{\fpcc}{\fpccpt}$ holds for a particular conjugacy class $\ccj$ if and only if the set $\mtc B_j$  contains exactly one element. This indicates that the corresponding character class relative to $\ccpt$ consists of a single element.
\er
Recall that $\cc$ is called a \emph{ $p$-fusion category} if $\fpcc=p^n$ for some $n\geq 1$.
\bc\label{p-nilp}
For any $p$-fusion category one has 
$$
\fp(\ccj)\mid \frac{\fpcc}{\fpccpt}.
$$
\ec
\bpf
Note that both $\fp(\ccj)$ and $\frac{\fpcc}{\fpccpt}$ are non-negative  integers that are powers of $p$-power. The previous inequality therefore implies the divisibility result.
\epf
\subsection{Proof of Theorem \ref{br-nil}} 
\bpf
Note that any nilpotent fusion category is weakly integral  because it can be viewed as iterated $G$-extensions of $\text{Vect}$. Furthermore, any braided nilpotent fusion category can be uniquely decomposed into a Deligne product of $p$-categories, as established in \cite[Theorem 1.1]{DGNO2}. Consequently, the divisibility property follows from Corollary \ref{p-nilp}. Indeed, if $\cc\simeq \prod_{p\mid\fpcc} \cc_p$ with each $\cc_p$ being a $p$-fusion category then it is straightforward to see  that $\czcc\simeq \prod_{p\mid\fpcc} \cz(\cc_p)$. Thus, each conjugacy class of $\cc$ is a Deligne product of a conjugacy class of each $\cc_p$. Moreover, $\ccpt\simeq \prod_{p\mid\fpcc} (\cc_p)_{pt}$ and the proof of the first part of the Theorem follows from Corollary \ref{p-nilp}.

Since $\cc$ is weakly integral Harada's identity applies and yields:
$$
\wtp^{\;2m-2}\big(\prod_{j=0}^mC_j \big)^2=\frac{ \wtp^{\;2m}\dim(\ccpt)}{\dimcc}\big(\sum_{j \in \mtc J_{\ccpt}} C_j\big)
$$
where $\wtp:=\prod_{j=0}^m \dim(\ccj)$. Given that the multiplicity of $C_0$ on the left hand side is an integer, it follows that $\dimcc/\dimccpt$ divides $\wtp^{\;2m}$. Additionally, since each $\dim(\ccj)$ divides this ratio, the proof is complete.
\epf
\subsection{Nilpotent modular categories}
Let $\cc$ be a modular category of Frobenius-Perron dimension $\fpcc=p^n$ for some prime number $p$. Since $\fp(\ccpt)$ divides $ \fpcc$, it follows that $\fp(\ccpt)=p^m$ for some $0\leq m\leq n$. 
 
Note that \cite[Lemma 2]{cps} establishes that in any modular tensor category of dimension $p^n$ one has that
$$
2 \leq m \leq n-3 
$$
Using Theorem \ref{gen-dual} the upper bound can be improved as follows:
\bl
Let $\cc$ be a an integral modular category of Frobenius-Perron dimension $\fpcc=p^n$ for some prime number $p$. Suppose that $\fp(\ccpt)=p^m$ for some $0\leq m\leq n$. Then either $2\leq m\leq n-4$ or the Frobenius-Perron dimension of the simple objects of $\cc$ are either $1$ or $p$.
\el
\bpf
Suppose that $m=n-3$. Note that $\fp(\ccj)=d_j^2$ and therefore, since $\fp(\ccj)\mid \frac{\fpcc}{\fp(\ccpt)} =p^{n-m}$ this implies that $d_j^2\mid p^3$ and therefore $d_j\in \{1,p\}$.

Additionally, consider the class equation:
$$
\fpcc=\sum_{j\in \mtcj_{\ccad}}\fp(\ccj)+\sum_{j\in \mtc J\setminus \mtcj_{\ccad}}\fp(\ccj).
$$
Since $p^2\mid \fp(\ccj)$ for all $j\notin \mtcj_{\ccad}$ it follows that $p^2\mid |\mtcj_{\ccad}|=|U(\cc)|$. On the other hand, since $\cc$ is modular, it follows that $|U(\cc)|=|\fp(\ccpt)|$, and therefore $m\geq 2$.
\epf
\bibliographystyle{alpha}
\bibliography{24nov}
\ed